\documentclass{etds}

\usepackage{amsfonts}
\usepackage{latexsym}
\usepackage{graphicx}

\newtheorem{thm}{Theorem}

\newtheorem{prop}[thm]{Proposition} 

\newtheorem{lemma}[thm]{Lemma}

\newcommand{\bprf}[1][Proof:]{\begin{list}{} 			{\setlength{\leftmargin}{1em} 			\setlength{\rightmargin}{0em}}                         \item {\bf \hspace{-1em}  #1 \ \ }} 
\begin{document}

\ETDS{1}{20}{30}{2003}
\runningheads{M. Pivato, R. Yassawi}{Limit Measures for Affine Cellular Automata II}

\title{Limit Measures for Affine Cellular Automata II}

\author{Marcus Pivato\affil{1}  and
 Reem Yassawi\affil{2}\footnote{This research partially supported by NSERC Canada.}}

\address{\affilnum{1}  Department of Mathematics, Trent University 
 \\ \email{{\tt pivato@xaravve.trentu.ca}} \\
\affilnum{2}  Department of Mathematics, Trent University,
 Lady Eaton College,\\
Peterborough, Ontario,  K9L 1Z6  Canada\\
\email{ {\tt ryassawi@trentu.ca}} }

\recd{???}

\begin{abstract}
  If ${\mathbb{M}}$ is a monoid, and ${\mathcal{ A}}$ is an abelian group, then ${\mathcal{ A}}^{\mathbb{M}}$ is
a compact abelian group; \ a {\em linear cellular automaton} (LCA) is
a continuous endomorphism ${\mathfrak{ F}}:{\mathcal{ A}}^{\mathbb{M}} {{\longrightarrow}} {\mathcal{ A}}^{\mathbb{M}}$ that commutes
with all shift maps.  If ${\mathfrak{ F}}$ is {\bf diffusive}, and $\mu$ is a {\bf
harmonically mixing} (HM) probability measure on ${\mathcal{ A}}^{\mathbb{M}}$, then the
sequence $\{{\mathfrak{ F}}^N \mu\}_{N=1}^{\infty}$ weak*-converges to the Haar measure
on ${\mathcal{ A}}^{\mathbb{M}}$, in density.  Fully supported Markov measures on
${\mathcal{ A}}^{\mathbb{Z}}$ are HM, and nontrivial LCA on ${\mathcal{ A}}^{\left({\mathbb{Z}}^D\right)}$ are
diffusive when ${\mathcal{ A}}={{\mathbb{Z}}_{/p}}$ is a prime cyclic group.

  In the present work, we provide sufficient conditions for diffusion
of LCA on ${\mathcal{ A}}^{\left({\mathbb{Z}}^D\right)}$ when ${\mathcal{ A}}={{\mathbb{Z}}_{/n}}$ is any cyclic
group or when ${\mathcal{ A}}=\left({{\mathbb{Z}}_{/p^r}}\right)^J$ ($p$ prime).  We also show that
any fully supported Markov random field on ${\mathcal{ A}}^{\left({\mathbb{Z}}^D\right)}$ is HM
(where ${\mathcal{ A}}$ is any abelian group).
\end{abstract}
\section{Introduction}

  Let $\displaystyle{\mathcal{ A}}$ be a finite abelian group, with discrete topology.  If
$\displaystyle{\mathbb{M}}$ is any set, then $\displaystyle{\mathcal{ A}}^{\mathbb{M}}$ is a compact abelian group when
endowed with the Tychonoff product topology and componentwise
addition.  If $\displaystyle{\mathbb{M}}$ is a monoid (for example, a lattice: \
$\displaystyle{\mathbb{Z}}^D \times {\mathbb{N}}^E$), then the action of $\displaystyle{\mathbb{M}}$ on itself by
translation induces a natural {\bf shift action} of $\displaystyle{\mathbb{M}}$ on
configuration space: \ for all $\displaystyle{\mathsf{ e}}\in{\mathbb{M}}$, and $\displaystyle{\mathbf{ a}}
\in {\mathcal{ A}}^{\mathbb{M}}$, define $\displaystyle {{{\boldsymbol{\sigma}}}^{{\mathsf{ e}}}} [{\mathbf{ a}}] \ = \ {\left[b_{\mathsf{ m}}  |_{{\mathsf{ m}}\in{\mathbb{M}}}^{} \right]}$
 where, $\displaystyle\forall {\mathsf{ m}}\in{\mathbb{M}}, \ \ \ b_{\mathsf{ m}} = a_{{\mathsf{ e}}.{\mathsf{ m}}}$.  Here ``$.$'' is the
monoid operator (``$+$'' for  $\displaystyle{\mathbb{M}} = {\mathbb{Z}}^D \times {\mathbb{N}}^E$).

  A {\bf linear cellular automaton} (LCA) is a continuous endomorphism
$\displaystyle{\mathfrak{ F}}:{\mathcal{ A}}^{\mathbb{M}} {\,\raisebox{0.3em}{$-$}\!\!\!\!\!\!\raisebox{-0.3em}{$\leftarrow$}\!\!\!\!\supset}$ which commutes with all shift maps.  If
$\displaystyle\mu$ is a measure on $\displaystyle{\mathcal{ A}}^{\mathbb{M}}$, it is natural to consider the
sequence of measures $\displaystyle{ \left\{{\mathfrak{ F}}^n \mu |_{_{{n\in{\mathbb{N}}}}} \right\} }$, and ask whether
this sequence converges in the weak* topology on the space
${{\mathcal{ M}}\left[{\mathcal{ A}}^{\mathbb{M}}\right] }$ of Borel probability measures on $\displaystyle{\mathcal{ A}}^{\mathbb{M}}$.  If
$\displaystyle{ \left\{{\mathfrak{ F}}^n \mu |_{_{{n\in{\mathbb{N}}}}} \right\} }$ does not itself converge, we might
hope at least for convergence in density (that is, convergence of a
subsequence $\displaystyle{ \left\{{\mathfrak{ F}}^j \mu |_{_{{j\in{\mathbb{J}}}}} \right\} }$, where $\displaystyle{\mathbb{J}} \subset {\mathbb{N}}$
is a subset of Ces\`aro  density $1$), or convergence of the
Ces\`aro  average $\displaystyle\frac{1}{N}\sum_{n=1}^N {\mathfrak{ F}}^n \mu$.

  Let $\displaystyle{\mathcal{ H}}$ denote the Haar measure on $\displaystyle{\mathcal{ A}}^{\mathbb{M}}$.  Since $\displaystyle{\mathcal{ H}}$
is invariant under the algebraic operations of $\displaystyle{\mathcal{ A}}^{\mathbb{M}}$, it seems
like a natural limit point for $\displaystyle{ \left\{{\mathfrak{ F}}^n \mu |_{_{{n\in{\mathbb{N}}}}} \right\} }$.  Indeed,
D. Lind showed \cite{Lind} that, if $\displaystyle{\mathcal{ A}} = {{\mathbb{Z}}_{/2}}$, and $\displaystyle{\mathfrak{ F}}$ is
the automaton defined: $\displaystyle{\mathfrak{ F}}({\mathbf{ a}})_0 = a_{(-1)} + a_{1}$, and $\displaystyle\mu$ is
any Bernoulli measure, then $\displaystyle {\mathbf{ w}}\!{\mathbf{ k}}^*\!\!-\!\!\!\lim_{N{\rightarrow}{\infty}}\frac{1}{N}\sum_{n=1}^N
{\mathfrak{ F}}^n \mu \ = \ {\mathcal{ H}}$.  Lind also showed that
${ \left\{{\mathfrak{ F}}^n \mu |_{_{{n\in{\mathbb{N}}}}} \right\} }$ does {\em not} converge to $\displaystyle{\mathcal{ H}}$;
convergence fails along the subsequence
 ${ \left\{{\mathfrak{ F}}^{\left(2^n\right)} \mu |_{_{{n\in{\mathbb{N}}}}} \right\} }$.

 Later, Ferrari, Maass, Martinez, and Ney showed similar
Ces\`aro  convergence results in a variety of special cases
\cite{FerMaassMartNey,MaassMartinezII}.
Recently, Pivato and Yassawi \cite{PivatoYassawi1} developed broad
sufficient conditions for convergence.  The concepts of {\bf harmonic
mixing} for measures and {\bf diffusion} for LCA were introduced;
\ if $\displaystyle\mu$ is a harmonically mixing probability measure
and $\displaystyle{\mathfrak{ F}}$ a diffusive LCA, then $\displaystyle{ \left\{{\mathfrak{ F}}^n \mu |_{_{{n\in{\mathbb{N}}}}} \right\} }$
weak* converges to $\displaystyle{\mathcal{ H}}$ in density, and thus, also in
Ces\`aro  mean.

\medskip

  This paper is a continuation of
\cite{PivatoYassawi1}.  First we will extend the
results on diffusion of LCA to a broader class of
abelian groups: \ in
\S\ref{S:diffusion.cyclic}, to the case when $\displaystyle{\mathcal{ A}} = {{\mathbb{Z}}_{/n}}$, for
any $\displaystyle n\in{\mathbb{N}}$, and then in \S\ref{S:diffusion.abelian}, to the
case when ${\mathcal{ A}} = \left({{\mathbb{Z}}_{/p^r}}\right)^J$ ($p$ prime, $J,r\in{\mathbb{N}}$).
Next, in
\S\ref{S:HM.MRF}, we demonstrate harmonic mixing for any Markov random
field on ${\mathcal{ A}}^{\left({\mathbb{Z}}^D\right)}$ with full support.  
\section{Preliminaries}

   We recommend that the reader consult
\cite{PivatoYassawi1} before reading the present work; \
we will depend heavily upon results introduced there.
We will now briefly review the relevant concepts;
\ all theorems in this section are proved in \cite{PivatoYassawi1}.

\subsection{Characters and Harmonic Mixing}

  Let ${{{\mathbb{T}}}^{1}}$ be the unit circle group.  A {\bf character} of
${\mathcal{ A}}^{\mathbb{M}}$ is a continuous group homomorphism $\phi:{\mathcal{ A}}^{\mathbb{M}} {{\longrightarrow}}
  {{{\mathbb{T}}}^{1}}$.  The set of all characters of ${\mathcal{ A}}^{\mathbb{M}}$ forms a group,
  denoted $\widehat{{\mathcal{ A}}^{\mathbb{M}}}$.
  
  If ${\left[\chi_{\mathsf{ m}}  |_{{\mathsf{ m}}\in{\mathbb{M}}}^{} \right]}$ is a sequence of characters of ${\mathcal{ A}}$,
with all but finitely many elements equal to the constant $1$-function
(denoted ``${{{\mathsf{ 1\!\!1}}}_{{}}}$''), then define $\displaystyle{\boldsymbol{\chi }} \ = \ \bigotimes_{{\mathsf{ m}}\in{\mathbb{M}}} \chi_{\mathsf{ m}}:{\mathcal{ A}}^{\mathbb{M}} {{\longrightarrow}}
{{{\mathbb{T}}}^{1}}$;\  thus, if ${\mathbf{ a}} \ = \ {\left[a_{\mathsf{ m}}  |_{{\mathsf{ m}}\in{\mathbb{M}}}^{} \right]}$ is an
element of ${\mathcal{ A}}^{\mathbb{M}}$, then \ $\displaystyle
{\boldsymbol{\chi }} ({\mathbf{ a}}) \ \ = \ \ \prod_{{\mathsf{ m}} \in {\mathbb{M}}}\chi_{\mathsf{ m}} (a_{\mathsf{ m}})$.
All elements of $\widehat{{\mathcal{ A}}^{\mathbb{M}}}$ arise in this manner.
The {\bf rank} of the character ${\boldsymbol{\chi }}$ is the number of
nontrivial entries in the {\bf coefficient system} ${\left[\chi_{\mathsf{ m}}  |_{{\mathsf{ m}}\in{\mathbb{M}}}^{} \right]}$.

  When ${\mathcal{ A}}={{\mathbb{Z}}_{/n}}$,  elements of $\widehat{{\mathcal{ A}}}$ are maps of the form
$\displaystyle \label{char.coef}
 \chi(a) \ = \ \exp\left(\frac{2 \pi {\mathbf{ i}}}{n} c\cdot a\right)$,
 where $c\in{{\mathbb{Z}}_{/n}}$ is some constant.
Elements of $\widehat{{\mathcal{ A}}^{\mathbb{M}}}$ are then products of the form
$\displaystyle{\boldsymbol{\chi }}({\mathbf{ a}}) \ = \ \prod_{{\mathsf{ m}}\in{\mathbb{M}}} \chi_{\mathsf{ m}} (a_{\mathsf{ m}})$, where,  \
$\forall {\mathsf{ m}}\in{\mathbb{M}}$, \ $\displaystyle \chi_{\mathsf{ m}}:a\mapsto \exp\left(\frac{2 \pi {\mathbf{ i}}}{n} c_{\mathsf{ m}}\cdot a\right)$ for some $c_{\mathsf{ m}}\in{\mathcal{ A}}$, with all but finitely many $c_{\mathsf{ m}}$
are equal to $0$.  In this case, we will use the term {\bf coefficient
system} also to describe the sequence ${\left[c_{\mathsf{ m}}  |_{{\mathsf{ m}}\in{\mathbb{M}}}^{} \right]}$.

  If $\mu$ is a measure on ${\mathcal{ A}}^{\mathbb{M}}$, then the {\bf Fourier
coefficients} of $\mu$ are defined: $\displaystyle {\widehat{\mu}}[{\boldsymbol{\chi }}] \ = \ {\left\langle {\boldsymbol{\chi }},\mu \right\rangle }
\ = \ \int_{{\mathcal{ A}}^{\mathbb{M}}} {\boldsymbol{\chi }} \ d \mu$, for every ${\boldsymbol{\chi }} \in
\widehat{{\mathcal{ A}}^{\mathbb{M}}}$.  The measure $\mu$ is called {\bf harmonically mixing}
if, for all $\epsilon > 0$, there is some $R > 0$ so that, for all ${\boldsymbol{\chi }}
\in \widehat{{\mathcal{ A}}^{\mathbb{M}}}$, \ \ 
$\displaystyle \left( \ \rule[-0.5em]{0em}{1em}       \begin{minipage}{40em}       \begin{tabbing}          ${{\sf rank}\left[{\boldsymbol{\chi }}\right]} \geq R$        \end{tabbing}      \end{minipage} \ \right) \
\ensuremath{\Longrightarrow} \
   \left( \ \rule[-0.5em]{0em}{1em}       \begin{minipage}{40em}       \begin{tabbing}          $\left| {\widehat{\mu}}[{\boldsymbol{\chi }}] \right| < \epsilon$        \end{tabbing}      \end{minipage} \ \right)$.

  Let ${{\mathcal{ M}}\left[{\mathcal{ A}}^{\mathbb{M}}; \ {\mathbb{C}}\right] }$ be the space of complex-valued Borel
measures on ${\mathcal{ A}}^{\mathbb{M}}$, treated as a Banach algebra under the total variation
norm, with operations of addition and convolution.  Let
${\mathcal{ H}} \subset {{\mathcal{ M}}\left[{\mathcal{ A}}^{\mathbb{M}}; \ {\mathbb{C}}\right] }$ be the set of harmonically mixing
measures.  

\begin{prop}{\sf \label{thm:mix.bernoulli.or.markov}}  
 Let ${\mathcal{ A}}$ be any finite abelian group.
\begin{enumerate}
\item ${\mathcal{ H}}$ is an ideal of ${{\mathcal{ M}}\left[{\mathcal{ A}}^{\mathbb{M}};\ {\mathbb{C}}\right] }$.

\item  ${\mathcal{ H}}$ is closed under the total variation norm and
dense in the weak* topology.

\item   ${\mathcal{ H}}$ contains all Bernoulli measures $\beta^{\otimes{\mathbb{M}}}$,
where $\beta$ is a measure on ${\mathcal{ A}}$ such that, for any subgroup ${\mathcal{ G}}
\subset {\mathcal{ A}}$, the support of $\mu$ extends over  more than one
coset of ${\mathcal{ G}}$.  

 \item If \ ${\mathbb{M}} = {\mathbb{Z}}$, then, for any $N>0$,  ${\mathcal{ H}}$ contains all $N$-step
Markov measures on ${\mathcal{ A}}^{\mathbb{Z}}$ giving nonzero probability to all
elements of ${\mathcal{ A}}^{\left[ 0..N \right]}$.

\item  ${\mathcal{ H}}$  contains any measure absolutely
continuous with respect to the aforementioned Bernoulli or Markov measures.
\hrulefill\ensuremath{\Box}
\end{enumerate}
 \end{prop}

\subsection{Linear Cellular Automata}

  A {\bf cellular automaton} (CA) is a continuous map
${\mathfrak{ F}}:{\mathcal{ A}}^{\mathbb{M}}{{\longrightarrow}}{\mathcal{ A}}^{\mathbb{M}}$ that commutes with all shift maps.  The
Curtis-Hedlund-Lyndon Theorem \cite{HedlundCA}
states that any CA is determined
by a {\bf local map} ${\mathfrak{ f}}:{\mathcal{ A}}^{\mathbb{U}}{{\longrightarrow}}{\mathcal{ A}}$, where ${\mathbb{U}}\subset{\mathbb{M}}$
is some finite subset (a ``neighbourhood of the identity'').
${\mathfrak{ F}}$ is an LCA if and only if ${\mathfrak{ f}}$ is
a homomorphism from the product group ${\mathcal{ A}}^{\mathbb{U}}$ into ${\mathcal{ A}}$.

 The set ${{{\mathbf{ E}}{\mathbf{ n}}{\mathbf{ d}}_{} \left({\mathcal{ A}}\right)}}$ of group endomorphisms from ${\mathcal{ A}}$ to itself is
a ring under composition and pointwise addition: \ if ${\mathfrak{ f}},{\mathfrak{ g}}$ are
in ${{{\mathbf{ E}}{\mathbf{ n}}{\mathbf{ d}}_{} \left({\mathcal{ A}}\right)}}$, then so are ${\mathfrak{ f}} \circ {\mathfrak{ g}}$ and ${\mathfrak{ f}} + {\mathfrak{ g}}$.
If ${\mathbf{Hom}_{} \left({\mathcal{ A}}^{\mathbb{U}}; \ {\mathcal{ A}}\right)}$ is the set of group homomorphisms from
${\mathcal{ A}}^{\mathbb{U}}$ into ${\mathcal{ A}}$, then there is a natural bijection between
$\left({{{\mathbf{ E}}{\mathbf{ n}}{\mathbf{ d}}_{} \left({\mathcal{ A}}\right)}}\right)^{\mathbb{U}}$ and ${\mathbf{Hom}_{} \left({\mathcal{ A}}^{\mathbb{U}}; \
{\mathcal{ A}}\right)}$, as follows: For each ${\mathsf{ u}}\in{\mathbb{U}}$, suppose that ${\mathfrak{ f}}_{\mathsf{ u}} \in
{{{\mathbf{ E}}{\mathbf{ n}}{\mathbf{ d}}_{} \left({\mathcal{ A}}\right)}}$.  Define ${\mathfrak{ f}}: {\mathcal{ A}}^{\mathbb{U}} {{\longrightarrow}} {\mathcal{ A}}$ by \ 
$\displaystyle
{\mathfrak{ f}} {\left[a_{\mathsf{ u}}  |_{{\mathsf{ u}}\in{\mathbb{U}}}^{} \right]}  \ \ = \ \  \sum_{{\mathsf{ u}}\in{\mathbb{U}}} {\mathfrak{ f}}_{\mathsf{ u}}(a_{\mathsf{ u}})$. \
Then ${\mathfrak{ f}}$ is a group homomorphism, and
every element of ${\mathbf{Hom}_{} \left({\mathcal{ A}}^{\mathbb{U}}; \ {\mathcal{ A}}\right)}$ arises in this manner.

  Thus, if ${\mathfrak{ F}}$ is an LCA, then there is some set of
{\bf coefficients} ${\left\{ {\mathfrak{ f}}_{\mathsf{ u}} \; ; \; {\mathsf{ u}}\in{\mathbb{U}} \right\} }$ so that,
for any ${\mathbf{ a}}\in{\mathcal{ A}}^{\mathbb{M}}$, ${\mathfrak{ F}}({\mathbf{ a}})={\mathbf{ b}}$, where
for all ${\mathsf{ m}}\in{\mathbb{M}}$, \ \ $\displaystyle
b_{\mathsf{ m}}  \ \ = \ \  \sum_{{\mathsf{ u}}\in{\mathbb{U}}} {\mathfrak{ f}}_{\mathsf{ u}}\left(a_{({\mathsf{ m}}.{\mathsf{ u}})}\right)
\ \ = \ \  \sum_{{\mathsf{ u}}\in{\mathbb{U}}} {\mathfrak{ f}}_{\mathsf{ u}}\left(  {{{\boldsymbol{\sigma}}}^{{\mathsf{ u}}}} ({\mathbf{ a}})_{\mathsf{ m}}\right)$.

For any ${\mathsf{ u}}\in{\mathbb{U}}$, treat ${\mathfrak{ f}}_{\mathsf{ u}}$ as an endomorphism on ${\mathcal{ A}}^{\mathbb{M}}$
by letting it act componentwise on elements of ${\mathcal{ A}}^{\mathbb{M}}$.  Then \
$\forall {\mathbf{ a}}\in{\mathcal{ A}}^{\mathbb{M}}$,\ \ $\displaystyle{\mathfrak{ F}}({\mathbf{ a}}) \ = \ \mbox{$\displaystyle \sum_{{\mathsf{ u}}\in{\mathbb{U}}}
{\mathfrak{ f}}_{\mathsf{ u}}\circ {{{\boldsymbol{\sigma}}}^{{\mathsf{ u}}}} ({\mathbf{ a}})$}$.  Thus, ${\mathfrak{ F}}$ can be written as a
formal ``polynomial of shift maps'':
\[
{\mathfrak{ F}} = \sum_{{\mathsf{ u}}\in{\mathbb{U}}} {\mathfrak{ f}}_{\mathsf{ u}} \circ  {{{\boldsymbol{\sigma}}}^{{\mathsf{ u}}}} .
\]
If ${\mathcal{ A}}={{\mathbb{Z}}_{/n}}$ ($n\in{\mathbb{N}}$), then the elements of
${{{\mathbf{ E}}{\mathbf{ n}}{\mathbf{ d}}_{} \left({\mathcal{ A}}\right)}}$ are all maps of the form ${\mathfrak{ f}}\left([a]_n\right)= [f\cdot
a]_n$, where ``$[\bullet]_n$'' refers to a mod-$n$ congruence
class, and $f \in {{\mathbb{Z}}_{/n}}$ is a constant, with multiplication via
the natural ring structure on ${{\mathbb{Z}}_{/n}}$.  In this case, we can
write \ $\displaystyle{\mathfrak{ F}} = \sum_{{\mathsf{ u}}\in{\mathbb{U}}} f_{\mathsf{ u}} \cdot  {{{\boldsymbol{\sigma}}}^{{\mathsf{ u}}}} $, \
  a polynomial with coefficients in ${{\mathbb{Z}}_{/n}}$.
For example, if ${\mathbb{M}}={\mathbb{Z}}$ and ${\mathfrak{ F}} =  {{{\boldsymbol{\sigma}}}^{-1}}  + 3\circ  {{{\boldsymbol{\sigma}}}^{1}} 
+ 5  {{{\boldsymbol{\sigma}}}^{2}} $, \
then this means that ${\mathfrak{ F}}({\mathbf{ a}})_k \ = \ \left[a_{(k-1)} + 3\cdot a_{(k+1)}
+ 5\cdot a_{(k+2)} \right]_n$.

\subsection{Diffusion}

  If ${\mathfrak{ F}}:{\mathcal{ A}}^{\mathbb{M}} {\,\raisebox{0.3em}{$-$}\!\!\!\!\!\!\raisebox{-0.3em}{$\leftarrow$}\!\!\!\!\supset}$ is an LCA and ${\boldsymbol{\chi }}$ is a character of
${\mathcal{ A}}^{\mathbb{M}}$, then ${\boldsymbol{\chi }} \circ {\mathfrak{ F}}$ is also a character.  ${\mathfrak{ F}}$ is called {\bf
diffusive}\footnote{In \cite{PivatoYassawi1}, this was called {\bf diffusion
in density}.  Since diffusion in density is the only kind we will encounter
in this paper, we have opted for more concise terminology.} if,
 for every nontrivial ${\boldsymbol{\chi }} \in \widehat{{\mathcal{ A}}^{\mathbb{M}}}$, \ 
there is some subset ${\mathbb{J}}_\chi\subset{\mathbb{N}}$ of density 1 so that
$\displaystyle
\lim_{{j {\rightarrow} {\infty}}\atop{j\in{\mathbb{J}}_\chi}}  {{\sf rank}\left[ {\boldsymbol{\chi }} \circ {\mathfrak{ F}}^j \right]}\ \ = \ \  {\infty}$.  We will abbreviate this to \ ``$\displaystyle  {{\sf rank}\left[ {\boldsymbol{\chi }} \circ {\mathfrak{ F}}^N \right]}
{ -\!\!\!-\!\!\!-\!\!\!-\!\!\!\!\!\!\!\!\!\!\!  ^{{\scriptscriptstyle {\rm dense}}}_{{\scriptscriptstyle N{\rightarrow}{\infty}}}   \!\!\!\!\!\!\!\!\!\longrightarrow } {\infty}$''.
\begin{thm}{\sf  \label{unboundedthm}}    Let $p$ be a prime number, and ${\mathcal{ A}} = {\mathbb{Z}}_{/p}$.  Let $D \geq 1$.  Then
  any nontrivial LCA on ${\mathcal{ A}}^{\left({\mathbb{Z}}^D\right)}$
  is diffusive.  \hrulefill\ensuremath{\Box} \end{thm}
\medskip

  By {\bf nontrivial} we mean that ${\mathfrak{ F}}$, as a polynomial of
shift maps, has more than one nontrivial coefficient.
The significance of diffusion and harmonic mixing is the 
following:

\begin{thm}{\sf \label{diffuse.and.mix.converge}}  
 Let ${\mathcal{ A}}$ be a finite abelian group, and ${\mathbb{M}}$ a countable monoid.
 Suppose that ${\mathfrak{ F}}:{\mathcal{ A}}^{\mathbb{M}} {\,\raisebox{0.3em}{$-$}\!\!\!\!\!\!\raisebox{-0.3em}{$\leftarrow$}\!\!\!\!\supset}$ is an LCA, and that
 $\mu$ is a harmonically mixing  measure on ${\mathcal{ A}}^{\mathbb{M}}$.

If ${\mathfrak{ F}}$ is diffusive, then there is
 a  set ${\mathbb{J}} \subset {\mathbb{N}}$ of Ces\`aro   density 1 
so that
  $\displaystyle {\mathbf{ w}}\!{\mathbf{ k}}^*\!\!-\!\!\!\lim_{{j {\rightarrow} {\infty}} \atop {j \in {\mathbb{J}}}} {\mathfrak{ F}}^j \mu
\  = \  {\mathcal{ H}}$.\

Thus,
$\displaystyle {\mathbf{ w}}\!{\mathbf{ k}}^*\!\!-\!\!\!\lim_{N{\rightarrow}{\infty}} \frac{1}{N} \sum_{n=1}^N {\mathfrak{ F}}^n \mu
\  =  \  {\mathcal{ H}}$. \hrulefill\ensuremath{\Box}
 \end{thm}

\medskip

  For example, ${\mathfrak{ F}}^n \mu$ weak*-converges to Haar measure in
density whenever $\mu$ is one of the aforementioned Bernoulli
or $N$-step Markov measures.  
\section{Diffusion on other cyclic groups
\label{S:diffusion.cyclic}}

   Suppose $n = p_1^{r_1} \cdot \ldots p_J^{r_J}$, where $p_1,\ldots,p_J$
are distinct primes and $r_1,\ldots,r_J \in {\mathbb{N}}$.   Let
${\mathcal{ A}} = {{\mathbb{Z}}_{/n}}$,  and ${\mathcal{ A}}_j = {{\mathbb{Z}}_{/q_j}}$,  with $q_j=p_j^{r_j}$, for
$j\in{\left[ 1..J \right]}$.  Then
$\displaystyle
{\mathcal{ A}} \ \cong \ \bigoplus_{j=1}^J {\mathcal{ A}}_j$, \
and thus,
$\displaystyle \widehat{{\mathcal{ A}}} \ \cong \
\bigoplus_{j=1}^J  \widehat{{\mathcal{ A}}_j}$. \ 
There is then a canonical identification: \ 
$\displaystyle {\mathcal{ A}}^{\mathbb{M}} \ \cong \ \bigoplus_{j=1}^J {\mathcal{ A}}_j^{\mathbb{M}}$, \
 and thus,\ $\displaystyle\widehat{{\mathcal{ A}}^{\mathbb{M}}}
\ \cong \ 
 \bigoplus_{j=1}^J \widehat{{\mathcal{ A}}_j^{\mathbb{M}}}$. \ \  Concretely: \ 
if ${\boldsymbol{\chi }}\in \widehat{{\mathcal{ A}}^{\mathbb{M}}}$ has coefficient sequence ${\left[c_{\mathsf{ m}}  |_{{\mathsf{ m}}\in{\mathbb{M}}}^{} \right]}$,
then, $\displaystyle {\boldsymbol{\chi }} \cong \bigoplus_{j=1}^J {\boldsymbol{\chi }}^{[j]}$, where
 for each $j\in{\left[ 1..J \right]}$, \ 
${\boldsymbol{\chi }}^{[j]}\in \widehat{{\mathcal{ A}}_j^{\mathbb{M}}}$ has coefficient sequence
${\left[c^{[j]}_{\mathsf{ m}}  |_{{\mathsf{ m}}\in{\mathbb{M}}}^{} \right]}$, with $c^{[j]}_{\mathsf{ m}} = [c_{\mathsf{ m}}]_{q_j}$
for all ${\mathsf{ m}}\in{\mathbb{M}}$.   Also,
\[
{{{\mathbf{ E}}{\mathbf{ n}}{\mathbf{ d}}_{} \left({\mathcal{ A}}\right)}}
\quad \cong \quad
\bigoplus_{i,j=1}^J {\mathbf{Hom}_{} \left({\mathcal{ A}}_i,{\mathcal{ A}}_j\right)} 
\quad  \raisebox{-1ex}{$\overline{\overline{{\scriptscriptstyle{\mathrm{(*)}}}}}$} \quad 
\bigoplus_{j=1}^J {{{\mathbf{ E}}{\mathbf{ n}}{\mathbf{ d}}_{} \left({\mathcal{ A}}_j\right)}}
\quad \cong \quad
\bigoplus_{j=1}^J {{\mathbb{Z}}_{/q_j}}. 
\]
 ($(*)$ is because cross-terms are trivial).
  Concretely, if $f \in {\mathbb{N}}$, and ${\mathfrak{ f}}:{\mathcal{ A}}{\,\raisebox{0.3em}{$-$}\!\!\!\!\!\!\raisebox{-0.3em}{$\leftarrow$}\!\!\!\!\supset}$ is
the map $[a]_n \mapsto [f\cdot a]_n$, then
${\mathfrak{ f}} = {\mathfrak{ f}}^{[1]} \oplus \ldots \oplus {\mathfrak{ f}}^{[J]}$, where, for each $j$,
${\mathfrak{ f}}^{[j]}:{\mathcal{ A}}_j{\,\raisebox{0.3em}{$-$}\!\!\!\!\!\!\raisebox{-0.3em}{$\leftarrow$}\!\!\!\!\supset}$ is
the map $[a]_{q_j} \mapsto [f_j \cdot a]_{q_j}$,
and where $f \equiv f_j \pmod p$. \ \
In particular, if $q_j$ divides $f$, then ${\mathfrak{ f}}^{[j]}$ is trivial.

  Thus, if ${\mathfrak{ F}}:{\mathcal{ A}}^{\mathbb{M}} {\,\raisebox{0.3em}{$-$}\!\!\!\!\!\!\raisebox{-0.3em}{$\leftarrow$}\!\!\!\!\supset}$ is the LCA 
$\displaystyle \sum_{{\mathsf{ u}}\in{\mathbb{U}}} {\mathfrak{ f}}_{\mathsf{ u}} \circ  {{{\boldsymbol{\sigma}}}^{{\mathsf{ u}}}} $,\
with ${\mathfrak{ f}}_{\mathsf{ u}} \in {{{\mathbf{ E}}{\mathbf{ n}}{\mathbf{ d}}_{} \left({\mathcal{ A}}\right)}}$, then, $\forall {\mathsf{ u}}\in{\mathbb{U}}$, we
can write $\displaystyle{\mathfrak{ f}}_{\mathsf{ u}} \ = \ {\mathfrak{ f}}_{\mathsf{ u}}^{[1]} \oplus \ldots \oplus {\mathfrak{ f}}_{\mathsf{ u}}^{[J]}$,
with ${\mathfrak{ f}}_{\mathsf{ u}}^{[j]} \in {{{\mathbf{ E}}{\mathbf{ n}}{\mathbf{ d}}_{} \left({\mathcal{ A}}_j\right)}}$ a scalar-multiplication map
determined by some  $f_{\mathsf{ u}}^{[j]} \in {{\mathbb{Z}}_{/q_j}}$,
and then write 
$\displaystyle {\mathfrak{ F}} =   \bigoplus_{j=1}^J {\mathfrak{ F}}^{[j]}$,\
where, $\forall j\in{\left[ 1..J \right]}$, \ 
${\mathfrak{ F}}^{[j]} : {\mathcal{ A}}_j^{\mathbb{M}} {\,\raisebox{0.3em}{$-$}\!\!\!\!\!\!\raisebox{-0.3em}{$\leftarrow$}\!\!\!\!\supset} $ \ 
is the LCA given by \ $\displaystyle
 \sum_{{\mathsf{ u}}\in{\mathbb{U}}} f^{[j]}_{\mathsf{ u}} \circ  {{{\boldsymbol{\sigma}}}^{{\mathsf{ u}}}} $. \
Note also that, if \
$\displaystyle {\boldsymbol{\chi }}  = \bigoplus_{j=1}^J {\boldsymbol{\chi }}_j \  \in \widehat{{\mathcal{ A}}^{\mathbb{M}}}$, \
then \
$\displaystyle {\boldsymbol{\chi }} \circ {\mathfrak{ F}}  =  \bigoplus_{j=1}^J\left( {\boldsymbol{\chi }}_j \circ {\mathfrak{ F}}^{[j]} \right)$.

\begin{lemma}{\sf \label{diffusive.direct.product}}  If \ $\displaystyle {\mathfrak{ F}}  = \bigoplus_{j=1}^J {\mathfrak{ F}}^{[j]}$ \ is an LCA
on  \ $\displaystyle \bigoplus_{j=1}^J {\mathcal{ A}}_j^{\mathbb{M}} $, \ then

\centerline{$\displaystyle
\left( \ \rule[-0.5em]{0em}{1em}       \begin{minipage}{40em}       \begin{tabbing}         ${\mathfrak{ F}}$ is diffusive         \end{tabbing}      \end{minipage} \ \right)  \iff 
 \left( \ \rule[-0.5em]{0em}{1em}       \begin{minipage}{40em}       \begin{tabbing}          $\forall j \in {\left[ 1..J \right]}$, \ \ 
${\mathfrak{ F}}^{[j]}$ is diffusive.         \end{tabbing}      \end{minipage} \ \right)$}
 \end{lemma}
\bprf

{\bf \hspace{-1em}  Proof of ``${\Longleftarrow}$'': \ \ }  Let ${\boldsymbol{\chi }} \in \widehat{{\mathcal{ A}}^{\mathbb{M}}}$ be nontrivial.  Thus,
${\boldsymbol{\chi }} = {\boldsymbol{\chi }}^{[1]} \oplus \ldots \oplus  {\boldsymbol{\chi }}^{[J]}$,
where at least one of ${\boldsymbol{\chi }}^{[j]} \in \widehat{{\mathcal{ A}}_j^{\mathbb{M}}}$,
is nontrivial;  suppose it is ${\boldsymbol{\chi }}^{[j_0]}$.  Since
${\mathfrak{ F}}^{[j_0]}$ is diffusive , we conclude: \
$\displaystyle
{{\sf rank}\left[{\boldsymbol{\chi }}\circ{\mathfrak{ F}}^N\right]} \geq  {{\sf rank}\left[{\boldsymbol{\chi }}^{[j_0]}\circ\left({\mathfrak{ F}}^{j_0}\right)^N\right]}
{ -\!\!\!-\!\!\!-\!\!\!-\!\!\!\!\!\!\!\!\!\!\!  ^{{\scriptscriptstyle {\rm dense}}}_{{\scriptscriptstyle n{\rightarrow}{\infty}}}   \!\!\!\!\!\!\!\!\!\longrightarrow }  {\infty}$.

{\bf \hspace{-1em}  Proof of ``$\ensuremath{\Longrightarrow}$'': \ \ }  Suppose that ${\mathfrak{ F}}_{j_0}$ is not diffusive.
Let ${\boldsymbol{\chi }}_{j_0}$ be some character on ${\mathcal{ A}}_{j_0}^{\mathbb{M}}$ so
that ${{\sf rank}\left[{\boldsymbol{\chi }}_{j_0}\circ{\mathfrak{ F}}_{j_0}\right]} \ensuremath{  -\!\!\!-\!\!\!-\hspace{-0.8em}\left/\rule[-0.5em]{0em}{1em}\right.\hspace{-1.4em}  ^{{\scriptscriptstyle {\rm dense}}}_{{\scriptscriptstyle n{\rightarrow}{\infty}}}   \!\!\!\!\!\!\!\!\!\longrightarrow }  {\infty}$, \
 and let $\displaystyle {\boldsymbol{\chi }} = \bigoplus_{j=1}^J \chi_j$,
where $\chi_j = {{{\mathsf{ 1\!\!1}}}_{{}}}$ for all $j\neq j_0$.  Then 
${{\sf rank}\left[{\boldsymbol{\chi }}\circ{\mathfrak{ F}}\right]} = {{\sf rank}\left[{\boldsymbol{\chi }}_{j_0}\circ{\mathfrak{ F}}_{j_0}\right]}
 \ensuremath{  -\!\!\!-\!\!\!-\hspace{-0.8em}\left/\rule[-0.5em]{0em}{1em}\right.\hspace{-1.4em}  ^{{\scriptscriptstyle {\rm dense}}}_{{\scriptscriptstyle n{\rightarrow}{\infty}}}   \!\!\!\!\!\!\!\!\!\longrightarrow }  {\infty}$, \
 so ${\mathfrak{ F}}$ is not diffusive.
 {\tt \hrulefill $\Box$ } \end{list}  \medskip

\medskip

  Hence, we have reduced the proof of diffusion to the prime power
case.

  Suppose ${\mathcal{ A}}={{\mathbb{Z}}_{/8}}$, and let ${\mathfrak{ F}} = {\mathbf{ Id}_{{}}}+2 {{{\boldsymbol{\sigma}}}^{1}} $
act on ${\mathcal{ A}}^{\mathbb{Z}}$.  Then ${\mathfrak{ F}}^{4\cdot N}={\mathbf{ Id}_{{}}}$ for all $N\in{\mathbb{N}}$,
so ${\mathfrak{ F}}$ cannot be diffusive.  This motivates
the conditions of the following theorem.

\begin{lemma}{\sf \label{diffusive.pr}}   Let ${\mathbb{M}} = {\mathbb{Z}}^D$.  Let ${\mathcal{ A}} = {{\mathbb{Z}}_{/q}}$, where $q=p^r$, \ $p$ is prime and $r\in{\mathbb{N}}$.  Let $\displaystyle
{\mathfrak{ F}} \ = \ \sum_{{\mathsf{ u}}\in{\mathbb{U}}} f_{\mathsf{ u}} \circ  {{{\boldsymbol{\sigma}}}^{{\mathsf{ u}}}} $. 

  If  $f_{\mathsf{ u}} \in {\left[ 0...q \right)}$ are relatively prime to $p$
for at least two ${\mathsf{ u}}\in{\mathbb{U}}$, then ${\mathfrak{ F}}$ is diffusive .
 \end{lemma}
\bprf
  Let ${\boldsymbol{\chi }} \in \widehat{{\mathcal{ A}}^{{\mathbb{M}}}}$ have
coefficient sequence ${\left[c_{\mathsf{ v}}  |_{{\mathsf{ v}}\in{\mathbb{V}}}^{} \right]}$, 
where 
$c_{\mathsf{ v}} \in {{\mathbb{Z}}_{/q}}$, for all ${\mathsf{ v}} \in {\mathbb{V}}$, with 
${\mathbb{V}}\subset {{\mathbb{M}}}$ some finite subset. Thus, ${\boldsymbol{\chi }}^{[N]} = {\boldsymbol{\chi }} \circ {\mathfrak{ F}}^N$ has
coefficient sequence ${\left[c_{\mathsf{ m}}^{[N]}  |_{{\mathsf{ m}}\in{{\mathbb{M}}}}^{} \right]}$, where,
for all $m\in{{\mathbb{M}}}$,
 \begin{equation} 
\label{coefficient.mod.pr}
 c_{\mathsf{ m}}^{[N]} \ = \ \sum_{{\mathsf{ v}}\in{\mathbb{V}}} \
	\sum_{{{\mathsf{ u}}_1,\ldots,{\mathsf{ u}}_N \in {\mathbb{U}}} \atop 
	{{\mathsf{ v}}+{\mathsf{ u}}_1+\ldots+{\mathsf{ u}}_N \ = \ {\mathsf{ m}}}}
 \ \ c_{\mathsf{ v}} \cdot f_{{\mathsf{ u}}_1} \cdot \ldots \cdot f_{{\mathsf{ u}}_N} \ \ 
 \end{equation} 
Thus, ${{\sf rank}\left[{\boldsymbol{\chi }} \circ {\mathfrak{ F}}^N\right]}$ is the
number of these coefficients that are nonzero, mod $q$.

{\bf Case 1:}
  {\em One of the coefficients ${ \left\{c_{\mathsf{ v}} |_{_{{{\mathsf{ v}}\in{\mathbb{V}}}}} \right\} }$ is
nonzero, mod $p$.}

 Consider the character ${\boldsymbol{\chi }}_{/p}$ and the (nontrivial)
LCA ${\mathfrak{ F}}_{/p}$ on ${{\mathbb{Z}}_{/p}}^{{\mathbb{M}}}$ induced by the
coefficients ${\left[c_{\mathsf{ v}}  |_{{\mathsf{ v}}\in{\mathbb{V}}}^{} \right]}$ and ${\left[f_{\mathsf{ u}}  |_{{\mathsf{ u}}\in{\mathbb{U}}}^{} \right]}$ respectively,
and, for all $N\in{\mathbb{N}}$, the character  ${\boldsymbol{\chi }}^{[N]}_{/p}$ induced by 
${\left[c^{[N]}_{\mathsf{ m}}  |_{{\mathsf{ m}}\in{\mathbb{M}}}^{} \right]}$.

  First, note that $\forall N\in{\mathbb{N}}$, \ ${\boldsymbol{\chi }}^{[N]}_{/p} =
{\boldsymbol{\chi }}_{/p} \circ {\mathfrak{ F}}^N_{/p}$ (simply consider equation
 (\ref{coefficient.mod.pr}), only mod $p$ instead).
  Notice that, for any $m$ and $N$,
if the expression in (\ref{coefficient.mod.pr})
is nonzero mod $p$, then it must be nonzero mod $q$.  Thus
${{\sf rank}\left[{\boldsymbol{\chi }}^{[N]}\right]} \geq {{\sf rank}\left[{\boldsymbol{\chi }}^{[N]}_{/p}\right]} = 
{{\sf rank}\left[{\boldsymbol{\chi }}_{/p} \circ {\mathfrak{ F}}_{/p}^N\right]}$.
Hence, it suffices to show that 
${{\sf rank}\left[{\boldsymbol{\chi }}_{/p} \circ {\mathfrak{ F}}_{/p}^N \right]} { -\!\!\!-\!\!\!-\!\!\!-\!\!\!\!\!\!\!\!\!\!\!  ^{{\scriptscriptstyle {\rm dense}}}_{{\scriptscriptstyle N{\rightarrow}{\infty}}}   \!\!\!\!\!\!\!\!\!\longrightarrow } {\infty}$.

 But one of ${ \left\{c_{\mathsf{ v}} |_{_{{{\mathsf{ v}}\in{\mathbb{V}}}}} \right\} }$ is nonzero, mod $p$, so ${\boldsymbol{\chi }}_{/p}$
is nontrivial as a character on ${{\mathbb{Z}}_{/p}}^{{\mathbb{M}}}$.  Thus, by Theorem
\ref{unboundedthm}, ${{\sf rank}\left[{\boldsymbol{\chi }}_{/p} \circ {\mathfrak{ F}}_{/p}^N\right]} { -\!\!\!-\!\!\!-\!\!\!-\!\!\!\!\!\!\!\!\!\!\!  ^{{\scriptscriptstyle {\rm dense}}}_{{\scriptscriptstyle N{\rightarrow}{\infty}}}   \!\!\!\!\!\!\!\!\!\longrightarrow } {\infty}$.

{\bf Case 2:}  {\em All the
  coefficients ${ \left\{c_{\mathsf{ v}} |_{_{{{\mathsf{ v}}\in{\mathbb{V}}}}} \right\} }$ are divisible by $p$.}

  Let $p^s$ be the greatest power of $p$ that divides all elements
of  ${ \left\{c_{\mathsf{ v}} |_{_{{{\mathsf{ v}}\in{\mathbb{V}}}}} \right\} }$; \ clearly $s < r$.
  Let ${\widetilde{r}}=r-s$ and ${\widetilde{q}}= p^{{\widetilde{r}}}$,
 and let ${\widetilde{\mathcal{ A}}} =   {{\mathbb{Z}}_{/{\widetilde{q}}}}$.  We will reduce the problem
to consideration of an LCA on ${{\mathbb{Z}}_{/{\widetilde{q}}}}$, and
then apply {\bf Case 1}.

   For all ${\mathsf{ v}} \in
  {\mathbb{V}}$,\ let ${\widetilde{c}}_{\mathsf{ v}} \ = \ c_{\mathsf{ v}}/p^s$, and let $\widetilde{\boldsymbol{\chi }} \in
  \widehat{{\widetilde{\mathcal{ A}}}^{{\mathbb{M}}}}$ be the corresponding character.  Let 
$\widetilde{\mathfrak{ F}}$ be the LCA on ${\widetilde{\mathcal{ A}}}^{\mathbb{M}}$ having the same coefficients
as ${\mathfrak{ F}}$;  thus, $\widetilde{\boldsymbol{\chi }}^{[N]} = \widetilde{\boldsymbol{\chi }} \circ \widetilde{\mathfrak{ F}}^N$ has
coefficient sequence ${\left[{\widetilde{c}}_{\mathsf{ m}}^{[N]}  |_{{\mathsf{ m}}\in{{\mathbb{M}}}}^{} \right]}$, where,
for all $m\in{{\mathbb{M}}}$, \ $\displaystyle
\label{coefficient.mod.ptlr}
 {\widetilde{c}}_{\mathsf{ m}}^{[N]} \ = \ \sum_{{\mathsf{ v}}\in{\mathbb{V}}} \
	\sum_{{{\mathsf{ u}}_1,\ldots,{\mathsf{ u}}_N \in {\mathbb{U}}} \atop 
	{{\mathsf{ v}}+{\mathsf{ u}}_1+\ldots+{\mathsf{ u}}_N \ = \ {\mathsf{ m}}}}
 \ \ {\widetilde{c}}_{\mathsf{ v}} \cdot f_{{\mathsf{ u}}_1} \cdot \ldots \cdot f_{{\mathsf{ u}}_N}$.

  Clearly, for all $N\in{\mathbb{N}}$ and ${\mathsf{ m}}\in{\mathbb{M}}$, $c_{\mathsf{ m}}^{[N]} \ = \
p^s\cdot {\widetilde{c}}_{\mathsf{ m}}^{[N]}$, so if ${\widetilde{c}}_{\mathsf{ m}}^{[N]} \not\equiv 0
\pmod{{\widetilde{q}}}$, then $c_{\mathsf{ m}}^{[N]} \not\equiv 0 \pmod{q}$.  Thus,
${{\sf rank}\left[{\boldsymbol{\chi }}^{[N]} \right]} \geq {{\sf rank}\left[\widetilde{\boldsymbol{\chi }}^{[N]} \right]}$.  But by
construction, at least one coefficient of $\widetilde{\boldsymbol{\chi }}$ is nonzero, mod
$p$.  Thus, by {\bf Case 1}, we have: ${{\sf rank}\left[\widetilde{\boldsymbol{\chi }} \circ \widetilde{\mathfrak{ F}}^N \right]}
{ -\!\!\!-\!\!\!-\!\!\!-\!\!\!\!\!\!\!\!\!\!\!  ^{{\scriptscriptstyle {\rm dense}}}_{{\scriptscriptstyle N{\rightarrow}{\infty}}}   \!\!\!\!\!\!\!\!\!\longrightarrow } {\infty}$.
 {\tt \hrulefill $\Box$ } \end{list}  \medskip

\begin{thm}{\sf  \label{unboundedthm.nonprime}  }   Let $n \in {\mathbb{N}}$, and ${\mathcal{ A}} = {\mathbb{Z}}_{/n}$.  Let $D \geq 1$, 
and let ${\mathfrak{ F}}:{\mathcal{ A}}^{\left({\mathbb{Z}}^D\right)}{\,\raisebox{0.3em}{$-$}\!\!\!\!\!\!\raisebox{-0.3em}{$\leftarrow$}\!\!\!\!\supset}$ be an LCA such that,
for each prime divisor $p$ of $n$,
at least two coefficients of ${\mathfrak{ F}}$ are relatively prime to $p$.
Then ${\mathfrak{ F}}$ is diffusive.   \end{thm}
\bprf
  Write $n = p_1^{r_1} \cdot \ldots p_J^{r_J}$,
${\mathcal{ A}} = {\mathcal{ A}}_1\oplus\ldots\oplus{\mathcal{ A}}_J$, ${\mathfrak{ F}} = {\mathfrak{ F}}_1\oplus\ldots\oplus{\mathfrak{ F}}_J$
 as before.  By Lemma \ref{diffusive.direct.product}, it suffices to
show that each of ${\mathfrak{ F}}_1,\ldots,{\mathfrak{ F}}_J$ is diffusive.  By Lemma
\ref{diffusive.pr} and the hypothesis, this is the case.
 {\tt \hrulefill $\Box$ } \end{list}  \medskip

\section{Diffusion on finite abelian groups
\label{S:diffusion.abelian}}

  Now suppose ${\mathcal{ A}}$ is an arbitrary finite abelian group.
Then ${\mathcal{ A}}$ has a canonical decomposition: \ $\displaystyle
{\mathcal{ A}} \   =  \ \bigoplus_{k=1}^K \ \bigoplus_{j=1}^{J_k} {\mathcal{ A}}_{(k,j)}$, 
with ${\mathcal{ A}}_{(k,j)} \ = \ {{\mathbb{Z}}_{/q_{(k,j)}}}$; \  $q_{(k,j)} = p_k^{r_{(k,j)}}$,
where
$p_1,\ldots,p_K$ are distinct primes with $r_{(k,1)}
, r_{(k,2)} , \ldots , r_{(k,J_k)}$ natural numbers for
each $k\in{\left[ 1..K \right]}$.  

 We will assume that ${\mathcal{ A}}$ is of the
special form where, for all $k \in {\left[ 1..K \right]}$,
$r_{k,1}=\ldots=r_{k,J_k} = r_k$.  In other words, \ 
$\displaystyle
{\mathcal{ A}} \  =  \ \bigoplus_{k=1}^K {\mathcal{ A}}_k$, \
with ${\mathcal{ A}}_k \ = \ \left({{\mathbb{Z}}_{/q_k}}\right)^{J_{k}}$, where
$p_1,\ldots,p_K$ are distinct primes, with $q_k = p_k^{r_k}$,
and $r_k,J_{k}\in {\mathbb{N}}$.
Thus, as before, \ 
$\displaystyle
{{{\mathbf{ E}}{\mathbf{ n}}{\mathbf{ d}}_{} \left({\mathcal{ A}}\right)}}
\ = \
 \bigoplus_{j,k=1}^K \ 
 {\mathbf{Hom}_{} \left({\mathcal{ A}}_{j},\ {\mathcal{ A}}_{k}\right)} 
\ = \
 \bigoplus_{k=1}^K \ {{{\mathbf{ E}}{\mathbf{ n}}{\mathbf{ d}}_{} \left({\mathcal{ A}}_k\right)}}$, 
(cross-terms are trivial), and ${\mathcal{ A}}^{\mathbb{M}} \ \cong \ {\mathcal{ A}}_1^{\mathbb{M}} \oplus
\ldots \oplus {\mathcal{ A}}_K^{\mathbb{M}}$, so we can write any LCA
 ${\mathfrak{ F}}:{\mathcal{ A}}^{\mathbb{M}} {\,\raisebox{0.3em}{$-$}\!\!\!\!\!\!\raisebox{-0.3em}{$\leftarrow$}\!\!\!\!\supset}$ as a direct sum \ ${\mathfrak{ F}} = {\mathfrak{ F}}_1 \oplus \ldots
\oplus {\mathfrak{ F}}_K$, \ where ${\mathfrak{ F}}_k:{\mathcal{ A}}_k^{\mathbb{M}}{\,\raisebox{0.3em}{$-$}\!\!\!\!\!\!\raisebox{-0.3em}{$\leftarrow$}\!\!\!\!\supset}$. 

 By Lemma
\ref{diffusive.direct.product}, to prove ${\mathfrak{ F}}$ is diffusive, it
suffices to show that each of ${\mathfrak{ F}}_1,\ldots,{\mathfrak{ F}}_K$ is diffusive.
 Hence, we will assume from now on that ${\mathcal{ A}} =
\left({{\mathbb{Z}}_{/q}}\right)^{J}$, where $p$ is prime, $q=p^r$, and $J \in
{\mathbb{N}}$.  Elements of ${\mathcal{ A}}$ are
thought of as $J$-tuples of  ${{\mathbb{Z}}_{/q}}$-elements.  ${\mathcal{ A}}$ is a
$J$-dimensional module over the commutative ring\footnote{If $r=1$
then $q=p$ is prime, ${{\mathbb{Z}}_{/q}}$ is a field, and ${\mathcal{ A}}$ is a
${{\mathbb{Z}}_{/q}}$-vector space.  It may be helpful to keep this case in
mind in what follows.} ${{\mathbb{Z}}_{/q}}$.  The endomorphisms of ${\mathcal{ A}}$ as
an abelian group are just the ${{\mathbb{Z}}_{/q}}$-linear endomorphisms of
this ${{\mathbb{Z}}_{/q}}$-module, and are described by $J\times J$ matrices of
elements in ${{\mathbb{Z}}_{/q}}$.
\begin{lemma}{\sf \label{character.vs.automorphism}}  
Let ${\mathcal{ A}} = \left({{\mathbb{Z}}_{/q}}\right)^{J}$, where $p$ is prime and $q=p^r$.
\begin{enumerate}
  \item Any $\chi \in \widehat{{\mathcal{ A}}}$ \ is of the form:
$\chi({\mathbf{ a}}) \ = \ \exp\left(\frac{2\pi{\mathbf{ i}}}{q} \cdot {\left\langle {\mathbf{ c}},{\mathbf{ a}} \right\rangle } \right)$,
where ${\mathbf{ c}} = (c_1,\ldots,c_J) \in ({{\mathbb{Z}}_{/q}})^J$, and
for any ${\mathbf{ a}} = (a_1,\ldots,a_J) \in ({{\mathbb{Z}}_{/q}})^J$,
we define ${\left\langle {\mathbf{ c}},{\mathbf{ a}} \right\rangle } = c_1 a_1+\ldots+c_J a_J$.  
Thus, $\chi$ is nontrivial if and only if ${\mathbf{ c}} \neq 0$.

  \item If \ ${\mathfrak{ f}}\in{{{\mathbf{ E}}{\mathbf{ n}}{\mathbf{ d}}_{} \left({\mathcal{ A}}\right)}}$ has matrix ${\mathbf{ F}}$ with adjoint
$\,^{\dagger}\!{{\mathbf{ F}}}$, then $\chi\circ{\mathfrak{ f}}$ is
the character \ ${\mathbf{ a}}\mapsto \exp\left(\frac{2\pi{\mathbf{ i}}}{q} \cdot {\left\langle {\mathbf{ c}}',{\mathbf{ a}} \right\rangle } \right)$,
where ${\mathbf{ c}}' = \,^{\dagger}\!{{\mathbf{ F}}}\cdot{\mathbf{ c}}$.  

  In particular,  $\chi\circ{\mathfrak{ f}}$ is nontrivial if and only if
${\mathbf{ c}}$ is not in $\ker[\,^{\dagger}\!{{\mathbf{ F}}}]$.

  \item Let \ ${\mathfrak{ f}} \in {{{\mathbf{ A}}{\mathbf{ u}}{\mathbf{ t}}_{} \left({\mathcal{ A}}\right)}}$.  If $\chi\in{\widehat{\mathcal{ A}}}$ is
nontrivial then $\chi\circ{\mathfrak{ f}}$ is also nontrivial.
\hrulefill\ensuremath{\Box}
\end{enumerate}
 \end{lemma}

\medskip

  Let ${\mathbb{V}}\subset{\mathbb{M}}$ be a subset not containing $0$.  If 
$\displaystyle{\mathfrak{ G}} = \sum_{{\mathsf{ m}}\in{\mathbb{M}}} {\mathfrak{ g}}_{\mathsf{ m}}  {{{\boldsymbol{\sigma}}}^{{\mathsf{ m}}}} $ is an LCA on ${\mathcal{ A}}^{\mathbb{M}}$,
then a subset ${\mathbb{W}}\subset{\mathbb{M}}$ is called {\bf ${\mathbb{V}}$-separating}
for ${\mathfrak{ G}}$ if,
for every $\displaystyle {\mathsf{ w}}\in{\mathbb{W}}$, \ ${\mathfrak{ g}}_{\mathsf{ w}} \in{{{\mathbf{ A}}{\mathbf{ u}}{\mathbf{ t}}_{} \left({\mathcal{ A}}\right)}}$, but for all
${\mathsf{ v}} \in {\mathbb{V}}$, \ \ ${\mathfrak{ g}}_{({\mathsf{ w}}-{\mathsf{ v}})} = 0$.  
Intuitively, ${\mathbb{W}}$ indexes a set of nontrivial (indeed, automorphic)
coefficients of ${\mathfrak{ G}}$, separated from one another by ${\mathbb{V}}$-shaped
``gaps''.  If ${\mathbb{U}}={\mathbb{V}}\sqcup\{0\}$, and $\displaystyle{\boldsymbol{\chi }}=\bigotimes_{{\mathsf{ u}}\in{\mathbb{U}}}
\chi_{\mathsf{ u}}$ is a character, then we will show that these gaps ensure
that $\left({\boldsymbol{\chi }}\circ{\mathfrak{ G}}\right)_{\mathsf{ w}}$ is nontrivial, for all ${\mathsf{ w}}\in{\mathbb{W}}$.
We will then construct ${\mathbb{V}}$-separating sets for ${\mathfrak{ G}} = {\mathfrak{ F}}^N$.  This
argument was already used implicitly to prove Theorem 15 in
\cite{PivatoYassawi1}.

\begin{prop}{\sf \label{diffusion.condition}}  
 Let ${\mathbb{M}}={\mathbb{Z}}^D$.  An LCA ${\mathfrak{ F}}:{\mathcal{ A}}^{\mathbb{M}}{\,\raisebox{0.3em}{$-$}\!\!\!\!\!\!\raisebox{-0.3em}{$\leftarrow$}\!\!\!\!\supset}$
is diffusive\ if,  for every finite subset ${\mathbb{V}} \subset {\mathbb{M}}$ not containing zero, and every $R\in{\mathbb{N}}$,
there is a set ${\mathbb{J}}_{({\mathbb{V}};R)}\subset{\mathbb{N}}$ of density 1
so that, for all $j\in{\mathbb{J}}_{({\mathbb{V}};R)}$ there is a ${\mathbb{V}}$-separating
set ${\mathbb{W}}_j\subset {\mathbb{M}}$ for ${\mathfrak{ F}}^j$ with ${\sf card}\left[{\mathbb{W}}_j\right] > R$.
 \end{prop}
\bprf Suppose  ${\mathfrak{ F}}$ is not diffusive; \ thus, there is some character
$\displaystyle {\boldsymbol{\chi }} = \prod_{{\mathsf{ u}}\in{\mathbb{U}}} \chi_{\mathsf{ u}}$ so that $\displaystyle {{\sf rank}\left[{\boldsymbol{\chi }}\circ{\mathfrak{ F}}^N\right]}
\ensuremath{  -\!\!\!-\!\!\!-\hspace{-0.8em}\left/\rule[-0.5em]{0em}{1em}\right.\hspace{-1.4em}  ^{{\scriptscriptstyle {\rm dense}}}_{{\scriptscriptstyle N{\rightarrow}{\infty}}}   \!\!\!\!\!\!\!\!\!\longrightarrow } {\infty}$; \ hence, there is some subset
${\mathbb{B}}\subset{\mathbb{N}}$ of nonzero upper density and some bound $R$
so that ${{\sf rank}\left[{\boldsymbol{\chi }}\circ{\mathfrak{ F}}^N\right]} < R$ for all $N\in{\mathbb{B}}$.

Fix ${\mathsf{ u}}_0\in{\mathbb{U}}$ and let ${\mathbb{V}} = {\left\{ {\mathsf{ u}}-{\mathsf{ u}}_0 \; ; \; {\mathsf{ u}}\in{\mathbb{U}}\setminus\{{\mathsf{ u}}_0\} \right\} }$;
let ${\mathbb{J}}_{({\mathbb{V}};R)}$ be the set described by the hypothesis. The set
${\mathbb{B}}\subset{\mathbb{N}}$ has nonzero upper density, so
${\mathbb{B}}\cap{\mathbb{J}}_{({\mathbb{V}};R)}\neq\emptyset$;\  let $j\in
{\mathbb{B}}\cap{\mathbb{J}}_{({\mathbb{V}};R)}$, and let ${\mathbb{W}}_j\subset{\mathbb{M}}$ be the
${\mathbb{V}}$-separating set for ${\mathfrak{ F}}^j$.

Write \ $\displaystyle {\mathfrak{ F}}^j \ = \ \sum_{{\mathsf{ m}}\in{\mathbb{M}}} {\mathfrak{ f}}^{[j]}_{\mathsf{ m}}  {{{\boldsymbol{\sigma}}}^{{\mathsf{ m}}}} $,
and then write \ 
$\displaystyle \chi\circ{\mathfrak{ F}}^j \ = \ \prod_{{\mathsf{ m}}\in{\mathbb{M}}} \chi^{[j]}_{\mathsf{ m}}$,\
\ where \
 $\displaystyle \chi^{[j]}_{\mathsf{ m}} \ = \  \prod_{{\mathsf{ u}}\in{\mathbb{U}}} \left(\chi_{\mathsf{ u}} \circ {\mathfrak{ f}}^{[j]}_{{\mathsf{ m}}-{\mathsf{ u}}}\right)$.
Then $\forall {\mathsf{ w}}\in{\mathbb{W}}_j$, \ \
$\displaystyle \chi^{[j]}_{({\mathsf{ w}}+{\mathsf{ u}}_0)}
 \ = \ \prod_{{\mathsf{ u}}\in{\mathbb{U}}} \left(\chi_{\mathsf{ u}} \circ {\mathfrak{ f}}^{[j]}_{({\mathsf{ w}}+{\mathsf{ u}}_0-{\mathsf{ u}})}\right)
\  = \ \left(\chi_{{\mathsf{ u}}_0} \circ {\mathfrak{ f}}^{[j]}_{{\mathsf{ w}}}\right)\cdot
 \prod_{{\mathsf{ v}}\in{\mathbb{V}}} \left(\chi_{({\mathsf{ v}}+{\mathsf{ u}}_0)} \circ {\mathfrak{ f}}^{[j]}_{({\mathsf{ w}}-{\mathsf{ v}})}\right) 
\ = \
\left(\chi_{{\mathsf{ u}}_0} \circ {\mathfrak{ f}}^{[j]}_{{\mathsf{ w}}}\right)\cdot \left(\prod_{{\mathsf{ v}}\in{\mathbb{V}}} {{{\mathsf{ 1\!\!1}}}_{{}}}\right)
\ = \
\left(\chi_{{\mathsf{ u}}_0} \circ {\mathfrak{ f}}^{[j]}_{{\mathsf{ w}}}\right)$, \
which is nontrivial by
 Lemma \ref{character.vs.automorphism}, because
${\mathfrak{ f}}^{[j]}_{{\mathsf{ w}}}$ is an automorphism.  Thus, $\chi^{[j]}_{{\mathsf{ w}}}
\neq {{{\mathsf{ 1\!\!1}}}_{{}}}$ for all ${\mathsf{ w}}\in{\mathbb{W}}+{\mathsf{ u}}_0$, a set of cardinality greater than
$R$, contradicting the hypothesis that  ${{\sf rank}\left[{\boldsymbol{\chi }}\circ{\mathfrak{ F}}^j\right]} < R$.
 {\tt \hrulefill $\Box$ } \end{list}  \medskip  

Applying Proposition \ref{diffusion.condition} often involves
tracking binomial coefficients, mod $p$, via Lucas' Theorem
\cite{PivatoYassawi1}.  For a fixed prime $p$, and any $n\in{\mathbb{N}}$,
let ${\mathbb{P}}(n)\in{\left[ 0...p \right)}^{\mathbb{N}}$ be the {\bf $p$-ary expansion} of
$n$ (conventionally written with digits in reversed order).  Thus, for
example, if $p=3$, then ${\mathbb{P}}(34) \ = \ \ldots0000001021$.

  If $n,N\in{\mathbb{N}}$, with ${\mathbb{P}}(n)= {\left[n^{[i]}  |_{i=0}^{{\infty}} \right]}$
and ${\mathbb{P}}(N)= {\left[N^{[i]}  |_{i=0}^{{\infty}} \right]}$
then we write ``$n\ll N$'' if $n^{[i]} \leq N^{[i]}$ for
all $i \in {\mathbb{N}}$.  Lucas' Theorem then implies:
\[
 \left( \ \left[ N \atop n \right]_p \not= \ 0 \ \right) \ 
\iff \left( \ \rule[-0.5em]{0em}{1em}       \begin{minipage}{40em}       \begin{tabbing}           $n \ll N $         \end{tabbing}      \end{minipage} \ \right)
\]

  A {\bf commuting automorphism linear cellular automaton}
is an LCA of the form\
$\displaystyle {\mathfrak{ F}} \ = \ \sum_{{\mathsf{ u}}\in{\mathbb{U}}} {\mathfrak{ f}}_{\mathsf{ u}} \circ  {{{\boldsymbol{\sigma}}}^{{\mathsf{ u}}}} $,\ where
${ \left\{{\mathfrak{ f}}_{\mathsf{ u}} |_{_{{{\mathsf{ u}}\in{\mathbb{U}}}}} \right\} } \subset {{{\mathbf{ A}}{\mathbf{ u}}{\mathbf{ t}}_{} \left({\mathcal{ A}}\right)}}$ is a commuting
collection of automorphisms of ${\mathcal{ A}}$.  For example:
\begin{itemize}
  \item ${ \left\{{\mathfrak{ f}}_{\mathsf{ u}} |_{_{{{\mathsf{ u}}\in{\mathbb{U}}}}} \right\} }$ are simultaneously diagonalizable
automorphisms.  In other words, there is some ${{\mathbb{Z}}_{/q}}$-basis
${\mathcal{ B}} = \{{\mathbf{ b}}_1,\ldots,{\mathbf{ b}}_J\}$ for ${\mathcal{ A}}$, so that the elements of
${\mathcal{ B}}$ are eigenvectors for every element of  ${ \left\{{\mathfrak{ f}}_{\mathsf{ u}} |_{_{{{\mathsf{ u}}\in{\mathbb{U}}}}} \right\} }$, and
all eigenvalues are relatively prime to $p$.

  \item  There is some  ${\mathfrak{ f}} \in {{{\mathbf{ A}}{\mathbf{ u}}{\mathbf{ t}}_{} \left({\mathcal{ A}}\right)}}$ so that $\forall {\mathsf{ u}}\in{\mathbb{U}}$, \ \ ${\mathfrak{ f}}_{\mathsf{ u}} = {\mathfrak{ f}}^{n_{\mathsf{ u}}}$ for some $n_{\mathsf{ u}} \in {\mathbb{Z}}$.
\end{itemize}
\begin{thm}{\sf \label{diffusive.calca}}  
 If ${\mathfrak{ G}}:{\mathcal{ A}}^{\left({\mathbb{Z}}^D\right)}{\,\raisebox{0.3em}{$-$}\!\!\!\!\!\!\raisebox{-0.3em}{$\leftarrow$}\!\!\!\!\supset}$ is a commuting automorphism
LCA with two or more nontrivial coefficients, then ${\mathfrak{ G}}$ is diffusive.
 \end{thm}
\bprf We will use Proposition \ref{diffusion.condition};\
the argument is basically identical to the proof of Theorem
15 in \cite{PivatoYassawi1}, so we will only sketch it here.

Suppose ${\mathfrak{ G}} \ = \ {\mathfrak{ g}}_0  {{{\boldsymbol{\sigma}}}^{{\mathsf{ n}}_0}}  + {\mathfrak{ g}}_1
 {{{\boldsymbol{\sigma}}}^{{\mathsf{ n}}_1}}  + \ldots + {\mathfrak{ g}}_U  {{{\boldsymbol{\sigma}}}^{{\mathsf{ n}}_U}} $, where
${\mathfrak{ g}}_0,\ {\mathfrak{ g}}_1,\ \ldots, \
{\mathfrak{ g}}_U \in {{{\mathbf{ A}}{\mathbf{ u}}{\mathbf{ t}}_{} \left({\mathcal{ G}}\right)}}$ commute, and where ${\mathsf{ n}}_0,\
{\mathsf{ n}}_1,\ \ldots, \ {\mathsf{ n}}_U \in {\mathbb{Z}}^D$.
We can rewrite: \ $ {\mathfrak{ G}}  \ = \ {\mathfrak{ g}}_0 \circ \left( {\mathfrak{ F}} \circ  {{{\boldsymbol{\sigma}}}^{{\mathsf{ n}}_0}}  \right),$  \ 
where:
\[
{\mathfrak{ F}}  \ = \ 
  {\mathbf{ Id}_{{}}} \ + \ {\mathfrak{ f}}_1  {{{\boldsymbol{\sigma}}}^{{\mathsf{ m}}_1}}  \left(  {\mathbf{ Id}_{{}}} \ + \ {\mathfrak{ f}}_2  {{{\boldsymbol{\sigma}}}^{{\mathsf{ m}}_2}} 
   \left[\rule[-0.5em]{0em}{1em} \ldots \left( {\mathbf{ Id}_{{}}} + 
 {\mathfrak{ f}}_{U-1}  {{{\boldsymbol{\sigma}}}^{{\mathsf{ m}}_{U-1}}}  \left[ {\mathbf{ Id}_{{}}} \ + \ {\mathfrak{ f}}_U  {{{\boldsymbol{\sigma}}}^{{\mathsf{ m}}_U}} \right] \right) \ldots  \right] \right),
\]
and, for all $u \in {\left[ 1..U \right]}$, \ \ ${\mathsf{ m}}_{\mathsf{ u}} \ = \ 
 {\mathsf{ n}}_{\mathsf{ u}}  -{\mathsf{ n}}_{u-1}$, and $ {\mathfrak{ f}}_{\mathsf{ u}} \  = \ {\mathfrak{ g}}_{u-1}^{-1} \circ {\mathfrak{ g}}_{\mathsf{ u}}$.  We can do this because ${\mathfrak{ g}}_0,\ {\mathfrak{ g}}_1,\ \ldots, \
{\mathfrak{ g}}_U$ are automorphisms, and thus, invertible. It suffices to
show that ${\mathfrak{ F}}$ is diffusive.

  Let $J\in{\mathbb{N}}$.  The coefficients of ${\mathfrak{ F}}$ commute, so we can
employ the Binomial Theorem ---and thus, Lucas' Theorem ---to compute
the coefficients of ${\mathfrak{ F}}^J$, mod $p$.
\[
\mbox{Let} \ \ {\mathcal{ L}}^U(J) \ = \ {\left\{ 
\left[k_1,k_2,\ldots,k_U\right] \in {\mathbb{N}}^U \; ; \;  k_U \ll k_{U-1} \ll \ldots
k_2 \ll k_1 \ll J  \right\} }.
\]
Then \ $\displaystyle {\mathfrak{ F}}^{J} 
  =  \sum_{{\mathsf{ n}} \in {\mathbb{Z}}^D} {\mathfrak{ f}}^{[J]}_{\mathsf{ n}} \circ  {{{\boldsymbol{\sigma}}}^{{\mathsf{ n}}}} $,\ 
where \ $\displaystyle {\mathfrak{ f}}^{[J]}_{\mathsf{ n}} \ = \ 
 \sum _{{{\mathbf{ k}} \in {\mathcal{ L}}^U(J)} \atop{\left(k_1 {\mathsf{ m}}_1 + \ldots + k_U {\mathsf{ m}}_U\right) = {\mathsf{ n}}}}
   {\mathfrak{ f}}^{[J]}_{({\mathbf{ k}})}$, \ 
and, for any ${\mathbf{ k}} \ = \ \left[k_1,k_2,\ldots,k_U\right] \in {\mathbb{N}}^U$, we define\
\[
  {\mathfrak{ f}}^{[J]}_{({\mathbf{ k}})}  :=    
{  \left[ { J } \atop { k_1 } \right]_{p} } {  \left[ { k_1 } \atop { k_2 } \right]_{p} } \ldots   {  \left[ { k_{U-1} } \atop { k_U } \right]_{p} } 
{\mathfrak{ f}}_1^{k_1} \circ {\mathfrak{ f}}_2^{k_2} \circ\ldots\circ {\mathfrak{ f}}_U^{k_U}.
\]
(See \cite{PivatoYassawi1} for details.)

  Fix a finite subset ${\mathbb{V}}\subset{\mathbb{Z}}^D$ not containing $0$, and let
$R>0$; \ we want to build a ${\mathbb{V}}$-separating set for ${\mathfrak{ F}}^J$ of
cardinality $R$.  To do this, note that there is some
$\Gamma\in{\mathbb{N}}$ such that, if $J\in{\mathbb{N}}$ and ${\mathbb{P}}(J)$ contains at
least $R$ ``gaps'' of size at least $\Gamma$ (ie. sequences of
$\Gamma$ successive zeros, delimited by nonzero entries), then we can
construct a set ${\mathbb{W}}_J
\subset {\mathbb{Z}}^D$ with ${\sf card}\left[{\mathbb{W}}_J\right] \geq R$, so that:
\begin{enumerate}
\item For every
${\mathsf{ w}}\in{\mathbb{W}}_J$, \ there is a unique ${\mathbf{ k}} \in {\mathcal{ L}}^U(J)$ so that
$\left(k_1 {\mathsf{ m}}_1 + \ldots + k_U {\mathsf{ m}}_U\right) = {\mathsf{ w}}$;\  thus
 ${\mathfrak{ f}}^{[J]}_{\mathsf{ w}} \ = \ {\mathfrak{ f}}^{[J]}_{({\mathbf{ k}})}\in{{{\mathbf{ A}}{\mathbf{ u}}{\mathbf{ t}}_{} \left({\mathcal{ A}}\right)}}$.

\item For all ${\mathsf{ v}} \in {\mathbb{V}}$, there are no ${\mathbf{ k}} \in {\mathcal{ L}}^U(J)$ with
$\left(k_1 {\mathsf{ m}}_1 + \ldots + k_U {\mathsf{ m}}_U\right) = {\mathsf{ w}}-{\mathsf{ v}}$; thus
 \ \ ${\mathfrak{ f}}^{[J]}_{({\mathsf{ w}}-{\mathsf{ v}})} = 0$.
\end{enumerate}
  Thus, ${\mathbb{W}}_J$ is ${\mathbb{V}}$-separating for ${\mathfrak{ F}}^J$.
  By Birkhoff's Ergodic Theorem, the set ${\mathbb{J}}_{(\Gamma;R)}$ of $J\in{\mathbb{N}}$
with $R$ such $\Gamma$-gaps is a set of
Ces\`aro  density one.  Thus, we satisfy the conditions of Proposition
\ref{diffusion.condition}.
 {\tt \hrulefill $\Box$ } \end{list}  \medskip  

  To apply Proposition \ref{diffusion.condition} it is clearly sufficient to
construct sets ${\mathbb{J}}_{({\mathbb{V}};R)}$ for some increasing sequence of numbers
$R_1,R_2,\ldots {\rightarrow}{\infty}$, along with a sequence ${\mathbb{V}}_1,{\mathbb{V}}_2,\ldots$
so that, for any finite ${\mathbb{V}}\subset{\mathbb{M}}$ we have ${\mathbb{V}}\subset {\mathbb{V}}_k+{\mathsf{ m}}$
for some ${\mathsf{ m}}\in{\mathbb{M}}$ and $k\in{\mathbb{N}}$.  Also, it suffices to prove
that the LCA ${\mathfrak{ F}}^K$ is diffusive\ for some power $K>0$: \ for any ${\boldsymbol{\chi }}\in
\widehat{{\mathcal{ A}}^{\mathbb{M}}}$, and any $k\in{\left[ 0...K \right)}$, ${\boldsymbol{\chi }}\circ{\mathfrak{ F}}^k$ is also a
character; \  if ${\mathfrak{ F}}^K$ is diffusive, then
$\displaystyle{{\sf rank}\left[{\boldsymbol{\chi }}\circ{\mathfrak{ F}}^k\circ{\mathfrak{ F}}^{n\cdot K}\right]}{ -\!\!\!-\!\!\!-\!\!\!-\!\!\!\!\!\!\!\!\!\!\!  ^{{\scriptscriptstyle {\rm dense}}}_{{\scriptscriptstyle n{\rightarrow}{\infty}}}   \!\!\!\!\!\!\!\!\!\longrightarrow }{\infty}$
for every $k\in{\left[ 0...K \right)}$, which in turn implies that
$\displaystyle{{\sf rank}\left[{\boldsymbol{\chi }}\circ{\mathfrak{ F}}^{n}\right]}{ -\!\!\!-\!\!\!-\!\!\!-\!\!\!\!\!\!\!\!\!\!\!  ^{{\scriptscriptstyle {\rm dense}}}_{{\scriptscriptstyle n{\rightarrow}{\infty}}}   \!\!\!\!\!\!\!\!\!\longrightarrow }{\infty}$.

\begin{figure}[htbp]
{\footnotesize
\hspace{-5em}\[
\begin{array}{rrccccccccccccccccccccccc}
j: 	&\!\! \ldots&\!\! * &\!\! 1 &\!\! 0 &\!\! * &\!\! \ldots&\!\! * &\!\! 1 &\!\! 0 &\!\! * &\!\! \ldots &\!\! * &\!\! 0 &\!\! q  &\!\! 1 &\!\! * &\!\! \ldots&\!\! * &\!\! * &\!\! * &\!\! \ldots &\!\! *&\!\! * &\!\! *\\
w: 	&\!\! \ldots&\!\! 0 &\!\! 0 &\!\! 0 &\!\! 0 &\!\! \ldots&\!\! 0 &\!\! 0 &\!\! 1 &\!\! 0 &\!\! \ldots &\!\! 0 &\!\! 0 &\!\! 1  &\!\! 0 &\!\! 0 &\!\! \ldots&\!\! 0 &\!\! 0 &\!\! 0 &\!\! \ldots &\!\! 0&\!\! 0 &\!\! 0\\
v: 	&\!\! \ldots&\!\! 0 &\!\! 0 &\!\! 0 &\!\! 0 &\!\! \ldots&\!\! 0 &\!\! 0 &\!\! 0 &\!\! 0 &\!\! \ldots &\!\! 0 &\!\! 0 &\!\! 0  &\!\! 0 &\!\! 0 &\!\! \ldots&\!\! 0 &\!\! 0 &\!\! * &\!\! \ldots &\!\! *&\!\! * &\!\! *\\
j+w:	&\!\! \ldots&\!\! * &\!\! 1 &\!\! 0 &\!\! * &\!\! \ldots&\!\! * &\!\! 1 &\!\! 1 &\!\! * &\!\! \ldots &\!\! * &\!\! 1 &\!\! 0  &\!\! 1 &\!\! * &\!\! \ldots&\!\! * &\!\! * &\!\! * &\!\! \ldots &\!\! *&\!\! * &\!\! *\\
j+w-1:	&\!\! \ldots&\!\! * &\!\! 1 &\!\! 0 &\!\! * &\!\! \ldots&\!\! * &\!\! 1 &\!\! 1 &\!\! * &\!\! \ldots &\!\! * &\!\! 1 &\!\! 0  &\!\! * &\!\! * &\!\! \ldots&\!\! * &\!\! * &\!\! * &\!\! \ldots &\!\! *&\!\! * &\!\! *\\
j+w-v:	&\!\! \ldots&\!\! * &\!\! 1 &\!\! 0 &\!\! * &\!\! \ldots&\!\! * &\!\! 1 &\!\! 1 &\!\! * &\!\! \ldots &\!\! * &\!\! 1 &\!\! 0  &\!\! * &\!\! * &\!\! \ldots&\!\! * &\!\! * &\!\! * &\!\! \ldots &\!\! *&\!\! * &\!\! *\\
j+w-v-1:&\!\! \ldots&\!\! * &\!\! 1 &\!\! 0 &\!\! * &\!\! \ldots&\!\! * &\!\! 1 &\!\! 1 &\!\! * &\!\! \ldots &\!\! * &\!\! 1 &\!\! 0  &\!\! * &\!\! * &\!\! \ldots&\!\! * &\!\! * &\!\! * &\!\! \ldots &\!\! *&\!\! * &\!\! *\\
2w: 	&\!\! \ldots&\!\! 0 &\!\! 0 &\!\! 0 &\!\! 0 &\!\! \ldots&\!\! 0 &\!\! 1 &\!\! 0 &\!\! 0 &\!\! \ldots &\!\! 0 &\!\! 1 &\!\! 0  &\!\! 0 &\!\! 0 &\!\! \ldots&\!\! 0 &\!\! 0 &\!\! 0 &\!\! \ldots &\!\! 0&\!\! 0 &\!\! 0\\
2v: 	&\!\! \ldots&\!\! 0 &\!\! 0 &\!\! 0 &\!\! 0 &\!\! \ldots&\!\! 0 &\!\! 0 &\!\! 0 &\!\! 0 &\!\! \ldots &\!\! 0 &\!\! 0 &\!\! 0  &\!\! 0 &\!\! 0 &\!\! \ldots&\!\! 0 &\!\! * &\!\! * &\!\! \ldots &\!\! *&\!\! * &\!\! *\\
2w-2v: 	&\!\! \ldots&\!\! 0 &\!\! 0 &\!\! 0 &\!\! 0 &\!\! \ldots&\!\! 0 &\!\! 1 &\!\! 0 &\!\! 0 &\!\! \ldots &\!\! 0 &\!\! 0 &\!\! q  &\!\! q &\!\! q &\!\! \ldots&\!\! q &\!\! * &\!\! * &\!\! \ldots &\!\! *&\!\! * &\!\! *\\
2w-2v-1:&\!\! \ldots&\!\! 0 &\!\! 0 &\!\! 0 &\!\! 0 &\!\! \ldots&\!\! 0 &\!\! 1 &\!\! 0 &\!\! 0 &\!\! \ldots &\!\! 0 &\!\! 0 &\!\! q  &\!\! q &\!\! q &\!\! \ldots&\!\! q &\!\! * &\!\! * &\!\! \ldots &\!\! *&\!\! * &\!\! *\\
        &\!\!       &\!\!   &\!\!   &\!\!\uparrow   &\!\!   &\!\!       &\!\!   &\!\!   &\!\!\uparrow   &\!\!   &\!\!        &\!\!   &\!\!   &\!\!    &\!\!\uparrow   &\!\!   &\!\!       &\!\!   &\!\!   &\!\! \uparrow&\!\!  &\!\! \uparrow&\!\! \uparrow &\!\! \uparrow\\
\mbox{Digit \#} &\!\! \lefteqn{\ldots\ldots}   &\!\!   &\!\!   &\!\! i_2   &\!\!   \lefteqn{\ldots\ldots\ldots}&\!\!      &\!\!   &\!\!   &\!\! i_1  &\!\!   
\lefteqn{\ldots\ldots\ldots}       &\!\!  &\!\!   &\!\!   &\!\!   &\!\! i_0  &\!\! \lefteqn{\ldots\ldots\ldots}&\!\!     &\!\!   &\!\!   &\!\! L_v&\!\! \ldots &\!\! 2&\!\! 1 &\!\! 0\\
\end{array}
\]
}
\caption{The $p$-ary expansions of $j$, $w$, $v$, etc. 
 Here, ``$*$'' represents any digit in in ${\left[ 0...q \right]}$.
\label{fig:mnm}}
\end{figure}
  
  Proposition \ref{diffusion.condition} can be applied even when the
coefficients of ${\mathfrak{ F}}$ do not commute.  For example:

        \refstepcounter{thm}        {\em }        \begin{list}{} 	{\setlength{\leftmargin}{1em} 	\setlength{\rightmargin}{0em}}         \item         {\bf Example  \thethm:} 
Let ${\mathcal{ A}} = \left({{\mathbb{Z}}_{/p}}\right)^2$, and let ${\mathfrak{ F}}:{\mathcal{ A}}^{\mathbb{Z}}{\,\raisebox{0.3em}{$-$}\!\!\!\!\!\!\raisebox{-0.3em}{$\leftarrow$}\!\!\!\!\supset}$
have local map ${\mathfrak{ f}}:{\mathcal{ A}}^{\{0,1\}}{{\longrightarrow}}{\mathcal{ A}}$ given: \ 
\[
{\mathfrak{ f}}\left(\left[{x_0}\atop{y_0}\right], \ \left[{x_1}\atop{y_1}\right] \right) 
\ = \ 
\left(y_0, \ x_0 + y_1\right)
\ = \ 
{\left[\begin{array}{ccccccccccccccccccccccccr} 0&1\\1&0 \end{array}\right]}\cdot{\left[\begin{array}{ccccccccccccccccccccccccr} x_0\\y_0 \end{array}\right]} \ + \
{\left[\begin{array}{ccccccccccccccccccccccccr} 0&0\\0&1 \end{array}\right]}\cdot{\left[\begin{array}{ccccccccccccccccccccccccr} x_1\\y_1 \end{array}\right]}
 \] 
This invertible LCA was studied in \cite{MaassMartinezII},
where it was shown to take fully supported Markov measures to Haar measure
in the weak* Ces\`aro  limit.  Proposition 3.1 of \cite{MaassMartinezII}
can be reformulated as:
\[
{\mathfrak{ F}}^N \ = \ \sum_{m=0}^N {\mathfrak{ f}}^{[N]}_m  {{{\boldsymbol{\sigma}}}^{m}} , 
\ \ \mbox{where} \ \ 
 {\mathfrak{ f}}^{[N]}_m = {\left[\begin{array}{ccccccccccccccccccccccccr}   \varphi^{(N-2)}_m & \varphi^{(N-1)}_m \\
			\varphi^{(N-1)}_m & \varphi^{(N)}_m  \end{array}\right]},
\]
\[
\mbox{with} \ \ 
  \varphi^{(N)}_m \ = \
 {\left\{ \begin{array}{rcl}                                  \left[ \left(\frac{N+m}{2}\right) \atop m \right]_p && 
\mbox{if} \ \  m\equiv N \pmod 2\\
  0 && \mbox{if} \ \  m \not\equiv N \pmod 2                                \end{array}  \right.  }.
\]
Thus, if $m\equiv (N-1) \pmod 2$, then the matrix ${\mathfrak{ f}}^{[N]}_m$ is
antidiagonal, and an automorphism iff $\varphi^{(N-1)}_m = \left[ \left(\frac{N-1+m}{2}\right) \atop
m \right]_p \neq 0$, which, by Lucas' Theorem, occurs only when $m \ll
\frac{N-1 +m}{2}$.  If $m\equiv N\equiv (N-2) \pmod 2$, then matrix
${\mathfrak{ f}}^{[N]}_m$ is diagonal, and an automorphism iff $m \ll
\frac{N+m}{2}$ and $m \ll \frac{N-2+m}{2}$.

  As noted earlier, it suffices to prove that ${\mathfrak{ F}}^2$ is
diffusive.  So, fix ${\mathbb{V}} = {\left( 0\ldots2V \right]}\subset {\mathbb{Z}}$ and
$R>0$;\ we will find a set ${\mathbb{J}}_{({\mathbb{V}};R)}$ and, for all
$j\in{\mathbb{J}}_{({\mathbb{V}};R)}$ some ${\mathbb{W}}_j \subset {\mathbb{Z}}$ with ${\sf card}\left[{\mathbb{W}}_j\right]>R$, so
that $2{\mathbb{W}}_j$ is ${\mathbb{V}}$-separating for ${\mathfrak{ F}}^{2j}$.  In other words,
 $\forall w\in{\mathbb{W}}_j$, \ ${\mathfrak{ f}}^{[2j]}_{2w} \in{{{\mathbf{ A}}{\mathbf{ u}}{\mathbf{ t}}_{} \left({\mathcal{ A}}\right)}}$, but
$\forall v \in {\mathbb{V}}$, \ \ ${\mathfrak{ f}}^{[2j]}_{(2w-v)} = 0$.  This is equivalent
to:

{\em $\forall w\in{\mathbb{W}}_j$, $\varphi^{(2j)}_{2w} \neq 0 \neq
\varphi^{(2j-2)}_{2w}$, but for all even $v=2u\in{\mathbb{V}}$, \ 
$\varphi^{(2j)}_{(2w-v)} = 0 =
\varphi^{(2j-2)}_{(2w-v)}$, and for all odd $v=2u+1\in{\mathbb{V}}$,
$\varphi^{(2j-1)}_{(2w-v)} = 0$.}  This, in turn, is equivalent to:

{\em For all $w\in{\mathbb{W}}_j$,
 \begin{equation} 
\label{pary1}
 2w \ll j+w\ \  \mbox{and}\  \ 2 w \ll j+w-1,
 \end{equation} 
but for all $u\in{\left( 0\ldots V \right]}$,
\begin{eqnarray}
\nonumber
 2w-2u &\not\ll& j+w-u, \hspace{2em}  2 w-2u \not\ll j+w-u-1, \\
\label{pary2}
 \mbox{and} \ \ 2w-2u-1 &\not\ll&  j+w-u-1.
\end{eqnarray}
}
So, let  $q=p-1$,\ $L_V = \lceil\log_p(V)\rceil+1$ and $L_R = \lceil\log_2(R)\rceil$, and let
${\mathbb{J}}_{({\mathbb{V}};R)}$
be the set of all $j\in{\mathbb{N}}$ such that ${\mathbb{P}}(j)$ contains the word ``$0q1$'' somewhere
after the first $L_V$ digits,
and contains at least $L_R$ separate instances of the word
``$10$'' after the ``$0q1$''.  By Birkhoff's Ergodic Theorem, 
${\mathbb{J}}_{({\mathbb{V}};R)}\subset{\mathbb{N}}$ has density 1.

  Suppose $j \in {\mathbb{J}}_{({\mathbb{V}};R)}$; and suppose that ``$0q1$'' occurs at position
$i_0 > L_v$, while ``$10$'' occurs at positions 
$i_{(L_R)} > \ldots > i_2> i_1$.
Let $w$ to be a number so that ${\mathbb{P}}(w)$ contains the word 
``$010$'' at $i_0$, and contains either ``$01$'' or ``$00$'' at each of
$i_1, i_2,\ldots,i_{(L_R)}$, with zeros everywhere else.  Clearly, we can
construct $2^{L_R} > R$ distinct numbers $w$ of this kind; \ let
${\mathbb{W}}_j$ be the set of all such numbers.

 For example, if $w$ has ``$01$'' at $i_1$ and ``$00$'' at $i_2$, and
$v\in{\left[ 0...V \right]}$, then the $p$-ary expansions of the relevant numbers
are depicted in Figure \ref{fig:mnm}.
By inspection, one can see that equations (\ref{pary1}) and  (\ref{pary2})
are satisfied.  Clearly, this will be true for any choice of $w\in{\mathbb{W}}_j$ and
$v\in{\mathbb{V}}$.
  \hrulefill			        \end{list}   			 
\section{Harmonic Mixing of Markov Random Fields
\label{S:HM.MRF}}

\paragraph*{Notation:}

  Suppose ${\mathbf{ a}} ={\mathcal{ A}}^{\mathbb{M}}$, with ${\mathbf{ a}}={\left[a_{\mathsf{ m}}  |_{{\mathsf{ m}}\in{\mathbb{M}}}^{} \right]}$.
If ${\mathbb{V}}\subset{\mathbb{M}}$, then ${\mathbf{ a}}\raisebox{-0.3em}{$\left|_{{\mathbb{V}}}\right.$} = {\left[a_{\mathsf{ v}}  |_{{\mathsf{ v}}\in{\mathbb{V}}}^{} \right]}
\in {\mathcal{ A}}^{\mathbb{V}}$.  This determines a continuous map
${\mathbf{ pr}_{{{\mathbb{V}}}}}:{\mathcal{ A}}^{\mathbb{M}}\ni {\mathbf{ a}} \mapsto{\mathbf{ a}}\raisebox{-0.3em}{$\left|_{{\mathbb{V}}}\right.$}\in {\mathcal{ A}}^{\mathbb{V}}$; \
if $\mu\in{{\mathcal{ M}}\left[{\mathcal{ A}}^{\mathbb{M}}\right] }$, then let ${\mathbf{ pr}_{{{\mathbb{V}}}}}^*\,(\mu)$ be the
${\mathbb{V}}$-marginal projection of $\mu$ (so that, for any ${\mathcal{ U}}\subset{\mathcal{ A}}^{\mathbb{V}}$,
${\mathbf{ pr}_{{{\mathbb{V}}}}}^*(\mu)\, [{\mathcal{ U}}] \ = \ \mu\left[ {\mathcal{ U}} \times {\mathcal{ A}}^{{\mathbb{M}}\setminus{\mathbb{V}}}\right]$).
 
  If ${\mathbf{ b}}\in{\mathcal{ A}}^{\mathbb{V}}$, then ${\left\langle {\mathbf{ b}} \right\rangle } = {\left\{ {\mathbf{ a}}\in{\mathcal{ A}}^{\mathbb{M}} \; ; \; {\mathbf{ a}}\raisebox{-0.3em}{$\left|_{{\mathbb{V}}}\right.$} =
{\mathbf{ b}} \right\} }$ is the associated {\bf cylinder set}, and, if $\mu
\in{{\mathcal{ M}}\left[{\mathcal{ A}}^{\mathbb{M}}\right] }$, then $\mu[{\mathbf{ b}}]$ is the measure of this cylinder
set.  If ${\mathbb{W}} \subset {\mathbb{M}}$ is disjoint from ${\mathbb{V}}$, and
${\mathbf{ c}}\in{\mathcal{ A}}^{\mathbb{W}}$, then ${\mathbf{ b}}\underline{\ }{\mathbf{ c}} \in {\mathcal{ A}}^{{\mathbb{V}}\sqcup{\mathbb{W}}}$ is defined so
that $\left({\mathbf{ b}}\underline{\ }{\mathbf{ c}}\right)\raisebox{-0.3em}{$\left|_{{\mathbb{V}}}\right.$} = {\mathbf{ b}}$ and
$\left({\mathbf{ b}}\underline{\ }{\mathbf{ c}}\right)\raisebox{-0.3em}{$\left|_{{\mathbb{W}}}\right.$}={\mathbf{ c}}$; \ thus, ${\left\langle {\mathbf{ b}}\underline{\ }{\mathbf{ c}} \right\rangle } =
{\left\langle {\mathbf{ b}} \right\rangle }\cap{\left\langle {\mathbf{ c}} \right\rangle }$.

  Let ${\mathcal{ B}}({\mathbb{V}})$ be the sigma-subalgebra of ${\mathcal{ A}}^{\mathbb{M}}$ generated by coordinates
in ${\mathbb{V}}$; \ if $\phi\in{\mathbf{ L}}^1\left({\mathcal{ A}}^{\mathbb{M}},\mu\right)$, let 
$\ensuremath{\mathbf{ E}_{{\mathbb{V}}}\left[ \phi \right]}\in{\mathbf{ L}}^1\left({\mathcal{ A}}^{\mathbb{V}},\mu\right)$ be the 
{\bf conditional expectation}
of $\phi$ given ${\mathcal{ B}}({\mathbb{V}})$, which we regard as a function on ${\mathcal{ A}}^{\mathbb{V}}$.
If ${\mathbf{ b}}\in{\mathcal{ A}}^{\mathbb{V}}$, then the {\bf conditional probability
measure} of $\mu$, {\bf given} ${\mathbf{ b}}$, is the unique measure
$\mu_{\mathbf{ b}}\in{{\mathcal{ M}}\left[{\mathcal{ A}}^{{\mathbb{M}}}\right] }$
such that ${\left\langle \phi,\mu_{\mathbf{ b}} \right\rangle } = \ensuremath{\mathbf{ E}_{{\mathbb{V}}}\left[ \phi \right]}({\mathbf{ b}})$
for every $\phi\in{\mathbf{ L}}^1\left({\mathcal{ A}}^{\mathbb{M}},\mu\right)$.
The map ${\mathcal{ A}}^{\mathbb{V}}\ni{\mathbf{ b}}\mapsto\mu_{\mathbf{ b}}\in{{\mathcal{ M}}\left[{\mathcal{ A}}^{{\mathbb{M}}}\right] }$ is measurable, and, if $\mu_{\mathbb{V}} =
{\mathbf{ pr}_{{{\mathbb{V}}}}}^*(\mu)$ \ is the marginal projection of $\mu$ onto ${\mathcal{ A}}^{\mathbb{V}}$, then
 $\mu$ has the {\bf disintegration}
 \cite{SchwartzDisintegrate,Furstenberg}: \ $ \displaystyle \mu \ = \
 \int_{{\mathcal{ A}}^{\mathbb{V}}} \mu_{\mathbf{ b}} \ d\mu_{\mathbb{V}}[{\mathbf{ b}}]$. \
Note that ${\mathbf{ pr}_{{{\mathbb{V}}}}}(\mu_{\mathbf{ b}})=\delta_{\mathbf{ b}}$, the point mass at ${\mathbf{ b}}$. \ 
If ${\mathbb{W}}\subset{\mathbb{M}}\setminus{\mathbb{V}}$ and ${\mathbf{ c}}\in{\mathcal{ A}}^{\mathbb{W}}$, 
we will sometimes write $\mu_{\mathbf{ b}}[{\mathbf{ c}}]$ as ``$\mu\left[{\mathbf{ c}}\ \left/\rule[-0.3em]{0em}{1em}\!\!\right/\ {\mathbf{ b}}\right]$'',
or, if ${\mathbf{ a}}\in{\mathcal{ A}}^{\mathbb{M}}$ is a $\mu$-random configuration, as \
``$\displaystyle\mu\left[\frac{{\mathbf{ a}}\raisebox{-0.3em}{$\left|_{\mathbb{W}}\right.$}={\mathbf{ c}}}{{\mathbf{ a}}\raisebox{-0.3em}{$\left|_{\mathbb{V}}\right.$} = {\mathbf{ b}}}\right]$''
\subsection{Markov Processes} 
 
  Let $({\mathbf{ X}},{\mathcal{ X}})$ be a measurable space, and let
$\mu\in{{\mathcal{ M}}\left[{\mathbf{ X}}^{\mathbb{Z}}\right] }$ be a probability measure.  Let ${\mathbb{U}} =
{\left[ 0...U \right)} \subset {\mathbb{Z}}$.  If $n\in{\mathbb{Z}}$ and ${\mathbf{ x}}\in
{\mathbf{ X}}^{({\mathbb{U}}+n)}$, then  ${\left\langle {\mathbf{ x}} \right\rangle } = {\left\{ {\mathbf{ y}}\in 
{\mathbf{ X}}^{\mathbb{Z}} \; ; \; {\mathbf{ y}}\raisebox{-0.3em}{$\left|_{({\mathbb{U}}+n)}\right.$} = {\mathbf{ x}} \right\} }$, and 
$\mu_{\mathbf{ x}}$ is the conditional probability measure of $\mu$, given
${\mathbf{ x}}$.

  $\mu$ is the {\bf path distribution} of a (${\mathbf{ X}}$-valued, $U$-step,
nonstationary) {\bf Markov process} if, for any $n\in{\mathbb{Z}}$ and
${\mathbf{ x}}\in {\mathbf{ X}}^{({\mathbb{U}}+n)}$, events occuring after time $n+U$ are
independent of those occuring before time $n$, relative to
$\mu_{\mathbf{ x}}$: \  for any ${\mathbb{V}}_{p} \subset
{\left( -{\infty}...n \right)}$, ${\mathbb{V}}_{f} \subset
{\left[ U+n...{\infty} \right)}$, and ${\mathbf{ y}}_{p} \in {\mathbf{ X}}^{{\mathbb{V}}_{p}}$ and 
${\mathbf{ y}}_{f} \in {\mathbf{ X}}^{{\mathbb{V}}_{f}}$, we have $\displaystyle \mu_{\mathbf{ x}} \left[{\mathbf{ y}}_{p}
\underline{\ } {\mathbf{ y}}_{f}\right] \ = \ \mu_{\mathbf{ x}}\left[{\mathbf{ y}}_{p}\right] \cdot
\mu_{\mathbf{ x}}\left[{\mathbf{ y}}_{f}\right]$.

  Any $U$-step Markov process is entirely described by its
$(U+1)$-dimensional marginals $\mu_{\left[ n...U+n \right]} = {\mathbf{ pr}_{{{\left[ n...U+n \right]}}}}^*
[\mu]$ for all $n\in{\mathbb{Z}}$, which are called the ($U$-step) {\bf
transition probabilities} of $\mu$.  If ${\mathbf{ X}}$ is finite, then $\displaystyle
{{\mathcal{ M}}\left[{\mathbf{ X}};\ {\mathbb{R}}\right] }\cong{\mathbb{R}}^{\mathbf{ X}}$; \ if $U=1$, then the transition
probabilities \ $\displaystyle \left\{\mu_{\{n,n+1\}}\right\}_{n\in{\mathbb{Z}}}$ can be
encoded by a sequence of {\bf transition probability matrices} \
${\left\{ {\mathbf{ Q}}^{(n)} \in {\mathbb{R}}^{{\mathbf{ X}}\times{\mathbf{ X}}} \; ; \; n\in{\mathbb{Z}} \right\} }$ and {\bf state
distributions} ${\left\{ \eta_n\in{\mathbb{R}}^{\mathbf{ X}} \; ; \; n\in{\mathbb{Z}} \right\} }$ so that, for any
$n\in{\mathbb{Z}}$, \ $\eta_{(n+1)} = {\mathbf{ Q}}^{(n)}\cdot\eta_{n}$, and, for any
$x_n,x_{(n+1)}\in{\mathbf{ X}}$, \ \ $\mu_{\{n,n+1\}}\left[x_n,x_{(n+1)}\right] \ = \
{\mathbf{ Q}}^{(n)}_{\left(x_{(n+1)};\ x_n\right)}\cdot\eta_n\left(x_n\right)$.

 If $\mu_{\left[ n...U+n \right]} = \mu_{\left[ 0...U \right]}$ for all $n\in{\mathbb{Z}}$, then
$\mu$ is {\bf stationary}.  If ${\mathbf{ X}}$ is finite and $U=1$, this means
there is some ${\mathbf{ Q}}\in{\mathbb{R}}^{{\mathbf{ X}}\times{\mathbf{ X}}}$ and $\eta \in {{\mathcal{ M}}\left[{\mathbf{ X}}\right] }$ (with
${\mathbf{ Q}}\cdot\eta = \eta$) so that ${\mathbf{ Q}}^{(n)} = {\mathbf{ Q}}$ and $\eta_n = \eta$
for all $n\in{\mathbb{Z}}$.  We call $\eta$ the {\bf stationary state
distribution}.

  If ${\mathfrak{ M}}\subset
{{\mathcal{ M}}\left[{\mathbf{ X}}^{\left[ 0...n \right]}\right] }$ is a finite family of transition probabilities,
we say $\mu$ is {\bf ${\mathfrak{ M}}$-semistationary} if $\mu_{\left[ n...U+n \right]}\in
 {\mathfrak{ M}}$ for all $n\in{\mathbb{Z}}$.  When ${\mathbf{ X}}$ is finite and $U=1$, this means
 that there are some finite families ${\mathfrak{ Q}}$ and ${\mathfrak{ H}}$ of transition
 probability matrices and state distributions, respectively, for $\mu$
 so that, for any $\eta\in{\mathfrak{ H}}$ and ${\mathbf{ Q}}\in{\mathfrak{ Q}}$, \ ${\mathbf{ Q}}\cdot\eta \in
 {\mathfrak{ H}}$; \ we say {\bf ${\mathfrak{ Q}}$-semistationary}.

  If $\mu$ is ${\mathfrak{ M}}$-semistationary, then $\mu$ has {\bf full support} if
every element of ${\mathfrak{ M}}$ has full support on ${\mathbf{ X}}^{\left[ 0...U \right]}$; as a
consequence, $\mu$ assigns nonzero probability to every finite
cylinder set.  If ${\mathbf{ X}}$ is finite and $U=1$, this means that every
entry of every transition probability matrix in ${\mathfrak{ Q}}$ is nonzero.

  If $\mu\in{{\mathcal{ M}}\left[{\mathbf{ X}}^{\mathbb{Z}}\right] }$ is a Markov process, $u,w\in{\mathbf{ X}}$, and
$n\in{\mathbb{Z}}$, then the {\bf sandwich measure}
$\,_n\mu_{u}^w\in{{\mathcal{ M}}\left[{\mathbf{ X}}\right] }$ is defined so that, if
${\mathbf{ x}}={\left[x_n  |_{n\in{\mathbb{Z}}}^{} \right]}$ is a $\mu$-random sequence, then for any
${\mathcal{ V}}\subset{\mathbf{ X}}$, \ $\displaystyle\,_n\mu_{u}^w({\mathcal{ V}}) \ = \
\mu\left[\frac{x_{n+1}\in{\mathcal{ V}}}{(x_n=u)\&(x_{n+2}=w)}\right]$.

\subsection{Exponential Harmonic Mixing}

 If ${\mathcal{ A}}$ is a finite abelian group,
and $\mu\in{{\mathcal{ M}}\left[{\mathcal{ A}}^{\mathbb{M}}; \ {\mathbb{C}}\right] }$, we will say $\mu$ is {\bf exponentially
harmonically mixing} with {\bf decay parameter} $\lambda>0$ (or
``$\lambda$-EHM'') if, for all ${\boldsymbol{\chi }}\in\widehat{{\mathcal{ A}}^{\mathbb{M}}}$ with
${{\sf rank}\left[{\boldsymbol{\chi }}\right]}\geq R$, we have $\left|{\left\langle {\boldsymbol{\chi }},\mu \right\rangle }\right| < e^{-\lambda\cdot
R}$.  It is straightforward to verify the following

\begin{lemma}{\sf \label{integral.ehm}}    Suppose $({\mathbf{ X}},\rho)$ is a probability space, and
 ${\mathbf{ X}}\ni{\mathbf{ x}}\mapsto\nu_{\mathbf{ x}}\in{{\mathcal{ M}}\left[{\mathcal{ A}}^{\mathbb{M}}; \ {\mathbb{C}}\right] }$ is a
measurable function so that $\nu_{\mathbf{ x}}$ is $\lambda$-EHM
for all ${\mathbf{ x}} \in{\mathbf{ X}}$.  If $\phi:{\mathbf{ X}}{{\longrightarrow}}{\mathbb{C}}$ is measurable
and ${\left\| \phi \right\|_{{{\infty}}} }   =1$, 
 then $\displaystyle \int_{\mathbf{ X}} \phi(x) \cdot \nu_{\mathbf{ x}} \ d\rho\left[{\mathbf{ x}}\right]$
is also $\lambda$-EHM.\hrulefill\ensuremath{\Box}
 \end{lemma}

\medskip

 If $\mu$ is a stationary, fully supported $U$-step
Markov measure on ${\mathcal{ A}}^{\mathbb{Z}}$, then $\mu$ is harmonically mixing (
{\bf Part 4} of Proposition \ref{thm:mix.bernoulli.or.markov} in this
paper, or Corollary 10 of \cite{PivatoYassawi1}).  The same method
easily generalizes to show:

\begin{prop}{\sf \label{markov.proc.exponential.HM}}   
Suppose ${\mathcal{ A}}$ is a finite abelian group, and that \ ${\mathfrak{ M}} \subset
{{\mathcal{ M}}\left[{\mathcal{ A}}^{\left[ 0...n \right]}\right] }$ is a finite family of fully supported
transition probabilities.

\begin{enumerate} 
\item  There is a constant $\lambda>0$  determined
by ${\mathfrak{ M}}$, so that, if $\mu$ is any ${\mathfrak{ M}}$-semistationary
Markov process on ${\mathcal{ A}}^{\mathbb{Z}}$, 
 then $\mu$ is $\lambda$-EHM.

\item  In particular, if $\mu$ is a $1$-step
${\mathfrak{ Q}}$-semistationary Markov process with full support, then
\ $ \displaystyle  -\lambda \ = \ 
\frac{1}{2}\cdot\sup_{{\xi,\chi\in{\widehat{\mathcal{ A}}}}\atop{\chi\neq{{{\mathsf{ 1\!\!1}}}_{{}}}}} \ \ \sup_{{\mathbf{ Q}},{\mathbf{ P}} \in {\mathfrak{ Q}}}
 \log {\left\| {\xi}_\bullet \cdot \,^{\dagger}\!{{\mathbf{ Q}}} \cdot {\chi}_\bullet \cdot \,^{\dagger}\!{{\mathbf{ P}}} \right\|_{{{\infty}}} }   $, \
  where ${\xi}_\bullet$ is the diagonal matrix with elements of
$\xi$ along the diagonal (so that, for any $\phi\in{\mathbb{C}}^{\mathcal{ A}}$,
${\xi}_\bullet\phi$ is the result of multiplying $\xi$ 
and $\phi$ componentwise), and where ${\left\| \bullet \right\|_{{{\infty}}} }   $ is
the uniform operator norm.
\end{enumerate}
 \end{prop}
\bprf {\em (Sketch)}
  Proposition 8 in \cite{PivatoYassawi1}
showed that a stationary $1$-step Markov matrix was harmonically
mixing;\  in fact, the proof showed that
$\left|{\left\langle {\boldsymbol{\chi }},\mu \right\rangle }\right| < e^{-\lambda R}$ for all ${\boldsymbol{\chi }}\in\widehat{{\mathcal{ A}}^{\mathbb{Z}}}$ with
${{\sf rank}\left[{\boldsymbol{\chi }}\right]} = R$, where  \
 $\\ \displaystyle -\lambda := \frac{1}{2}\cdot\sup_{{\xi,\chi\in{\widehat{\mathcal{ A}}}}\atop{\chi\neq{{{\mathsf{ 1\!\!1}}}_{{}}}}}  
 \log {\left\| {\xi}_\bullet \cdot \,^{\dagger}\!{{\mathbf{ Q}}} \cdot {\chi}_\bullet \cdot \,^{\dagger}\!{{\mathbf{ Q}}} \right\|_{{{\infty}}} }   $.  The
same argument works for a semistationary $1$-step process;  this yields
{\bf Part 2}.

  The proof of Corollary 10 in \cite{PivatoYassawi1} showed how any
fully supported $U$-step process could be ``recoded'' as a fully
supported $1$-step process; \ harmonic mixing of the latter implied
harmonic mixing of the former.  Corollary 10 thus followed from
Proposition 8.  By an identical argument
{\bf Part 1} follows from {\bf Part 2}.
 {\tt \hrulefill $\Box$ } \end{list}  \medskip  

\subsection{Markov Random Fields
\label{S:mark.rand.field}}

\begin{figure}
\centerline{\includegraphics{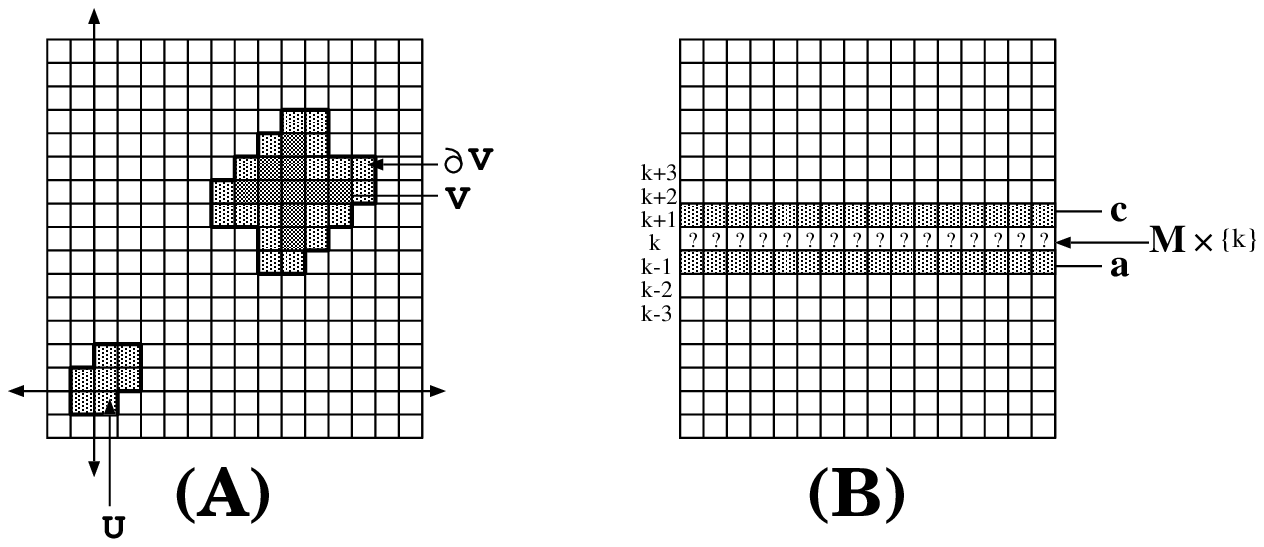}}
\caption{{\bf(A)} \ ${\mathbb{U}}$, ${\mathbb{V}}$, and $\partial {\mathbb{V}}$.
\quad {\bf(B)} \ A sandwich measure.
\label{fig:nhood}}
\end{figure}

 Let  ${\mathbb{U}}\subset{\mathbb{M}}$ be a finite ``neighbourhood of
$0$''  (e.g.
${\mathbb{M}}={\mathbb{Z}}^D$ and ${\mathbb{U}}={\left[ -1...1 \right]}^D$).  For any subset
${\mathbb{V}} \subset{\mathbb{M}}$, let $cl({\mathbb{V}}) := {\mathbb{V}} + {\mathbb{U}}$, and let $\partial({\mathbb{V}}) :=
cl({\mathbb{V}})\setminus {\mathbb{V}}$ (see Figure \ref{fig:nhood}{\bf(A)}). 

  $\mu \in{{\mathcal{ M}}\left[{\mathcal{ A}}^{\mathbb{M}}\right] }$ is a {\bf (nonstationary) Markov random
field} \cite{Bremaud,KindermannSnell} with {\bf interaction range} ${\mathbb{U}}$ 
(or ``${\mathbb{U}}$-MRF'') if, for any ${\mathbb{W}}\subset{\mathbb{M}}$,
and any ${\mathbf{ a}}\in{\mathcal{ A}}^{\partial({\mathbb{W}})}$, events occuring ``inside'' ${\mathbb{W}}$
are independent of those occuring ``outside'', relative to the conditional
measure $\mu_{\mathbf{ a}}$.  In other words, for any ${\mathbb{V}}_{in}\subset{\mathbb{W}}$,
${\mathbb{V}}_{out}
\subset {\mathbb{M}}\setminus cl({\mathbb{W}})$, and ${\mathbf{ b}}_{in} \in {\mathcal{ A}}^{{\mathbb{V}}_{in}}$,
 ${\mathbf{ b}}_{out} \in {\mathcal{ A}}^{{\mathbb{V}}_{out}}$, we have: \ 
$\displaystyle \mu_{\mathbf{ a}}\left[ {\mathbf{ b}}_{in} \underline{\ } {\mathbf{ b}}_{out} \right] \ = \
\mu_{\mathbf{ a}}\left[ {\mathbf{ b}}_{in}\right] \cdot \mu_{\mathbf{ a}}\left[ {\mathbf{ b}}_{out} \right]$.

  For example, if ${\mathbb{M}}={\mathbb{Z}}$, then the $U$-step Markov processes
on ${\mathcal{ A}}^{\mathbb{M}}$ are exactly the Markov random fields
with interaction range ${\mathbb{U}} = {\left( -U...U \right)}$.

  $\mu$ is {\bf stationary} if it is invariant under translation by
  ${\mathbb{M}}$.  In this case, $\mu_{({\mathbb{U}}+{\mathsf{ m}})} = \mu_{\mathbb{U}}$ for
every ${\mathsf{ m}}\in{\mathbb{M}}$, and $\mu_{\mathbb{U}} = {\mathbf{ pr}_{{{\mathbb{U}}}}}^*(\mu)$ is called
the {\bf local interaction} for $\mu$.

If ${\mathfrak{ I}} \subset {{\mathcal{ M}}\left[{\mathcal{ A}}^{\mathbb{U}}\right] }$ is finite, then
$\mu$ is {\bf ${\mathfrak{ I}}$-semistationary} if
$\mu_{({\mathbb{U}}+{\mathsf{ m}})} \in {\mathfrak{ I}}$ for every ${\mathsf{ m}}\in{\mathbb{M}}$.
${\mathfrak{ I}}$ is called the set of {\bf local interactions}.
We say $\mu$ has {\bf full support} if all elements of ${\mathfrak{ I}}$
have full support on ${\mathcal{ A}}^{\mathbb{U}}$.

\paragraph*{Lamination Processes:}
 
 Suppose ${\widetilde{\mathbb{U}}} \subset {\widetilde{\mathbb{M}}} = {\mathbb{M}}\times{\mathbb{Z}}$ and $\mu\in{{\mathcal{ M}}\left[{\mathcal{ A}}^{\widetilde{\mathbb{M}}}\right] }$ is
a ${\widetilde{\mathbb{U}}}$-MRF.  By a suitable recoding, we can assume ${\widetilde{\mathbb{U}}} = {\mathbb{U}} \times
\{-1,0,1\}$ for some ${\mathbb{U}}\subset{\mathbb{M}}$.  We can then realize $\mu$ via
an ${\mathcal{ A}}^{\mathbb{M}}$-valued, $1$-step Markov process, called the {\bf
lamination process}.  Intuitively, we imagine this Markov process as
constructing a $\mu$-random configuration in ${\mathcal{ A}}^{\widetilde{\mathbb{M}}}$ by laying
down successive random ``${\mathbb{M}}$-layers'', with each ${\mathbb{M}}$-layer
conditional on the previous one.

\begin{figure}
\centerline{\includegraphics{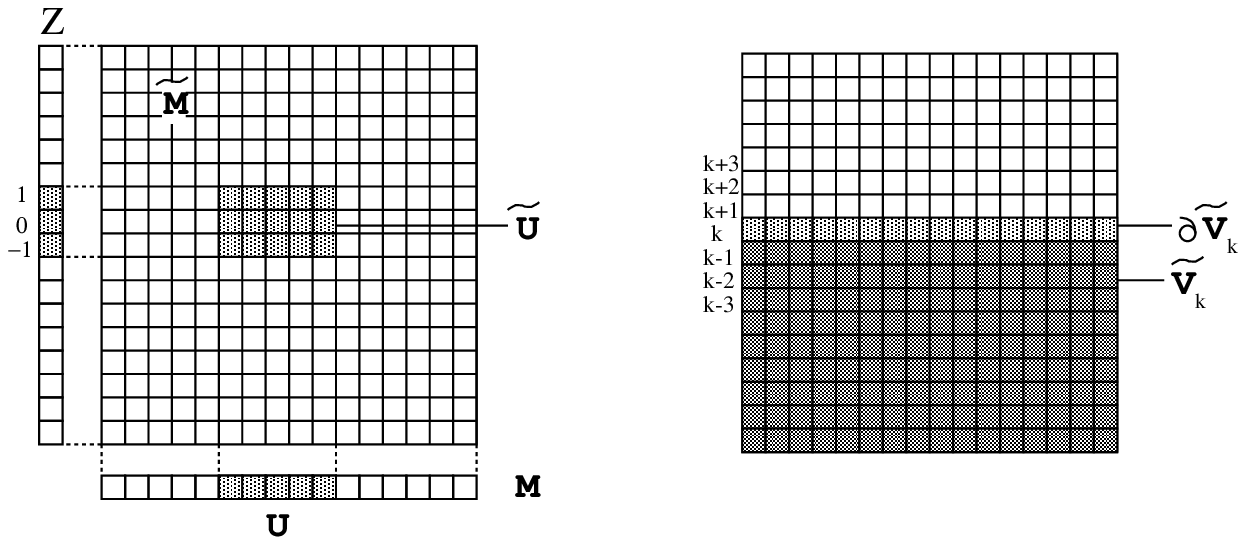}}
\caption{$ {\widetilde{\mathbb{M}}} = {\mathbb{M}}\times{\mathbb{Z}}$ and ${\widetilde{\mathbb{U}}} = {\mathbb{U}} \times \{-1,0,1\}$; \ \ 
 ${\widetilde{\mathbb{V}}}_k$ and $\partial {\widetilde{\mathbb{V}}}_k$.\label{fig:laminate}}
\end{figure}

   To see that this is a Markov process on ${\mathcal{ A}}^{\mathbb{M}}$, fix $k$,
and let ${\widetilde{\mathbb{V}}}_k \ = \ {\mathbb{M}} \times {\left( -{\infty}...k \right)}$ (the ``past'').  Then
$\partial({\widetilde{\mathbb{V}}}_k) = {\mathbb{M}} \times \{k\}$ (the ``present'') and ${\widetilde{\mathbb{M}}}
\setminus cl({\widetilde{\mathbb{V}}}_k) = {\mathbb{M}} \times
{\left( k...{\infty} \right)}$ (the ``future''); \  the Markov field condition of $\mu$ implies
that events in the past are independent of those in the future, given
complete information about the present (see Figure
\ref{fig:laminate}).
The original field measure $\mu\in{{\mathcal{ M}}\left[{\mathcal{ A}}^{{\mathbb{M}}\times{\mathbb{Z}}}\right] }$ is also
the path distribution (as a measure on $\left({\mathcal{ A}}^{\mathbb{M}}\right)^{\mathbb{Z}}$) for the
lamination process.  

\paragraph*{Sandwich Measures:}

  Again assume ${\widetilde{\mathbb{M}}} = {\mathbb{M}}\times{\mathbb{Z}}$ and ${\widetilde{\mathbb{U}}} = {\mathbb{U}} \times \{-1,0,1\}$.
If ${\mathbf{ a}}\in{\mathcal{ A}}^{{\mathbb{M}}\times\{k-1\}}$ and ${\mathbf{ c}}\in{\mathcal{ A}}^{{\mathbb{M}}\times\{k+1\}}$ (see
Figure \ref{fig:nhood}{\bf(B)}), then the {\bf sandwich measure} determined
by ${\mathbf{ a}}$ and ${\mathbf{ c}}$ is the sandwich measure $\,_{_{(k-1)}}\mu_{\mathbf{ a}}^{\mathbf{ c}}
\in{{\mathcal{ M}}\left[{\mathcal{ A}}^{\mathbb{M}}\right] }$ of
the lamination process; \ since $k$ is implicit in the definition of
${\mathbf{ a}}$ and ${\mathbf{ c}}$, we will suppress it, and denote the sandwich measure
as ``$\mu_{\mathbf{ a}}^{\mathbf{ c}}$''.  In other words, $\mu_{\mathbf{ a}}^{\mathbf{ c}}$ is the conditional measure
$\mu_{{\mathbf{ a}}\underline{\ }{\mathbf{ c}}}$, projected onto ${\mathcal{ A}}^{{\mathbb{M}}\times\{k\}}$.  The
following is easy to verify:

\begin{lemma}{\sf \label{sandwich.lemma}}  \begin{enumerate}

\item  $\mu_{\mathbf{ a}}^{\mathbf{ c}}$ is a Markov random field on ${\mathcal{ A}}^{\mathbb{M}}$, with 
interaction range ${\mathbb{U}}$.

\item If $\widetilde{\mathfrak{ I}}\subset{{\mathcal{ M}}\left[{\mathcal{ A}}^{\widetilde{\mathbb{U}}}\right] }$ and
 $\mu$ is $\widetilde{\mathfrak{ I}}$-semistationary, then there is some finite
${\mathfrak{ I}}\subset{{\mathcal{ M}}\left[{\mathcal{ A}}^{\mathbb{U}}\right] }$ so that all sandwich measures of $\mu$ are
${\mathfrak{ I}}$-semistationary.

\item If $\mu$ has full support, then so does every sandwich measure
of $\mu$.\hrulefill\ensuremath{\Box}
\end{enumerate}
 \end{lemma}

\medskip

  The harmonic mixing of an MRF depends on the the
harmonic mixing of its sandwich measures:

\begin{prop}{\sf \label{sandwich.unto.mixing}}  
 If ${\mathcal{ A}}$ is a finite abelian group, and
$\mu$ is a semistationary MRF on ${\mathcal{ A}}^{{\mathbb{M}}\times{\mathbb{Z}}}$
and all sandwich measures of $\mu$
are $\lambda$-EHM, then $\mu$ is $\lambda'$-EHM, where
$\lambda' = \lambda/2$.
 \end{prop}
\bprf See \S\ref{S:uhm}.
 {\tt \hrulefill $\Box$ } \end{list}  \medskip  

  From this follows our main result:

\begin{thm}{\sf }  
 Suppose ${\mathcal{ A}}$ is a finite abelian group, ${\mathbb{U}}\subset{\mathbb{Z}}^D$, and
 let $\widetilde{\mathfrak{ I}}\subset{{\mathcal{ M}}\left[{\mathcal{ A}}^{\mathbb{U}}\right] }$ be a finite set of local interactions
with full support.  Then $\exists\lambda>0$ so that
 if $\mu$ is any $\widetilde{\mathfrak{ I}}$-semistationary
 MRF on ${\mathcal{ A}}^{{\mathbb{Z}}^D}$, then $\mu$ is $\lambda$-EHM.
  \end{thm}
\bprf ({\em by induction on $D$}) 
\ \   If $D=1$, this is just Proposition \ref{markov.proc.exponential.HM}.

  Suppose inductively that the claim is true for MRFs on
${\mathbb{Z}}^{D-1}$, and let $\mu\in{\mathcal{ A}}^{{\mathbb{Z}}^D}$.  By Lemma
\ref{sandwich.lemma}, all sandwich measures of $\mu$ are
${\mathfrak{ I}}$-semistationary MRFs on ${\mathcal{ A}}^{{\mathbb{Z}}^{D-1}}$, where ${\mathfrak{ I}}$ is some
finite set of local interactions with full support.  Thus, by
induction hypothesis, all these sandwich measures are $\lambda$-EHM for
some $\lambda>0$.  Thus, by Proposition
\ref{sandwich.unto.mixing}, $\mu$ is $\lambda'$-EHM,
with $\lambda'=-\lambda/2$.
 {\tt \hrulefill $\Box$ } \end{list}  \medskip

\subsection{Markov Operators}

  When ${\mathbf{ X}}$ is finite, a $1$-step ${\mathbf{ X}}$-valued Markov process 
can be defined by a series of with transition probability matrices
$\{{\mathbf{ Q}}^{(n)}\}_{n\in{\mathbb{Z}}}$.  These matrices define linear operators
on the space ${{\mathcal{ M}}\left[{\mathbf{ X}}; \ {\mathbb{R}}\right] } \cong {\mathbb{R}}^{\mathbf{ X}}$, so that, if
$\eta_n\in{{\mathcal{ M}}\left[{\mathbf{ X}}\right] }$ is the state distribution at time $n$,
then ${\mathbf{ Q}}^{(n)}\cdot\eta_n=\eta_{n+1}$ is the state distribution at time $n+1$.

  When ${\mathbf{ X}}$ is an arbitrary measurable space (with sigma-algebra
${\mathcal{ X}}$), transition probabilities are described by linear operators on
the vector space ${{\mathcal{ M}}\left[{\mathbf{ X}};{\mathbb{R}}\right] }$ (which, for technical reasons,
we will treat as linear operators on ${{\mathcal{ M}}\left[{\mathbf{ X}}; \
{\mathbb{C}}\right] }$).

       \begin{list}{}
	{\setlength{\leftmargin}{1em}
	\setlength{\rightmargin}{0em}}
        \item {\bf Idea:}
 Informally speaking, a {\bf Markov operator} is linear operator
${\mathbf{ Q}}:{{\mathcal{ M}}\left[{\mathbf{ X}};{\mathbb{C}}\right] }{\,\raisebox{0.3em}{$-$}\!\!\!\!\!\!\raisebox{-0.3em}{$\leftarrow$}\!\!\!\!\supset}$ mapping the set ${{\mathcal{ M}}\left[{\mathbf{ X}}\right] }$ of
probability measures into itself.  Suppose $(y_0,y_1)\in{\mathbf{ X}}^{\{0,1\}}$
is a random couple, and ${\mathbf{ Q}}$ is the transition probability operator
from time $0$ to time $1$.   If $x\in{\mathbf{ X}}$, and $\delta_x\in{{\mathcal{ M}}\left[{\mathbf{ X}}\right] }$ is
the point mass at $x$, then  the probability measure ${\mathbf{ q}}_x :=
{\mathbf{ Q}}(\delta_x)$  is the conditional state distribution of $y_1$
 given that $y_0=x$: \ for all ${\mathcal{ U}}\subset{\mathbf{ X}}$, \ $\displaystyle {\mathbf{ q}}_x[{\mathcal{ U}}] \ = \
\ensuremath{{\sf Prob}\left[ \frac{y_1\in{\mathcal{ U}}}{y_0=x} \right] }$.  When ${\mathbf{ X}}$ is finite, measures on
${\mathbf{ X}}$ are vectors and ${\mathbf{ Q}}$ is a matrix, and ${\mathbf{ q}}_x$ is just the $x$th
column of this matrix.

  Suppose $y_0,y_1$ have distributions $\eta_0,\eta_1 \in{{\mathcal{ M}}\left[{\mathbf{ X}}\right] }$
respectively, with $\eta_{1}={\mathbf{ Q}}(\eta_0)$.  If
$\phi:{\mathbf{ X}}{{\longrightarrow}}{\mathbb{C}}$ is a measurable function, then the expected value
of $\phi(y_{1})$ is given by ${\left\langle \phi, \ \eta_{1} \right\rangle } \ = \ {\left\langle \phi,
{\mathbf{ Q}}(\eta_0) \right\rangle }
\ = \ {\left\langle \,^{\dagger}\!{{\mathbf{ Q}}}(\phi), \ \eta_0 \right\rangle }$, where $\,^{\dagger}\!{{\mathbf{ Q}}}$
 is the {\bf adjoint} of ${\mathbf{ Q}}$.  

  For any measurable ${\mathcal{ U}}\subset{\mathbf{ X}}$, let $\,^{\dagger}\!{{\mathbf{ q}}}\rule[-0.5em]{0em}{1em}_{\mathcal{ U}} := \,^{\dagger}\!{{\mathbf{ Q}}}\left({{{\mathsf{ 1\!\!1}}}_{{{\mathcal{ U}}}}}\right)$.
Thus, for any $x\in{\mathbf{ X}}$, \ $\,^{\dagger}\!{{\mathbf{ q}}}\rule[-0.5em]{0em}{1em}_{\mathcal{ U}}(x)  = 
\,^{\dagger}\!{{\mathbf{ Q}}}\left({{{\mathsf{ 1\!\!1}}}_{{{\mathcal{ U}}}}}\right)\left(x\right)  =  
{\left\langle \,^{\dagger}\!{{\mathbf{ Q}}}\left({{{\mathsf{ 1\!\!1}}}_{{{\mathcal{ U}}}}}\right), \ \delta_x \right\rangle }  =  
{\left\langle {{{\mathsf{ 1\!\!1}}}_{{{\mathcal{ U}}}}}, \ {\mathbf{ Q}}(\delta_x) \right\rangle }  =  
{\left\langle {{{\mathsf{ 1\!\!1}}}_{{{\mathcal{ U}}}}}, \ {\mathbf{ q}}_x \right\rangle }  =  {\mathbf{ q}}_x[{\mathcal{ U}}]$.  When ${\mathbf{ X}}$ is finite and ${\mathbf{ Q}}$ is a matrix and ${\mathcal{ U}} = \{u\}$ is a singleton set, then $\,^{\dagger}\!{{\mathbf{ q}}}\rule[-0.5em]{0em}{1em}_{\mathcal{ U}}$ is just the
$u$th row of ${\mathbf{ Q}}$ (or the $u$th column of $\,^{\dagger}\!{{\mathbf{ Q}}}$).
\end{list}  

We  need to develop some technology to make these ideas well-defined.
\paragraph*{Formalism:}
If $\Phi:{\mathbf{ X}}{{\longrightarrow}}{\mathbb{C}}$ is measurable, then let \ 
$\displaystyle {\left\| \Phi \right\|_{{{\infty}}} }    =
\sup_{{\mathbf{ x}}\in{\mathbf{ X}}} |\Phi(x)|$, \  and consider the Banach space \
$\displaystyle  {\mathcal{ M}}_{\infty}({\mathbf{ X}},{\mathcal{ X}}) = \
\{\Phi:{\mathbf{ X}}{{\longrightarrow}}{\mathbb{C}} \ ; \ \Phi\ \mbox{ measurable,} \\ {\left\| \Phi \right\|_{{{\infty}}} }   <{\infty}\}$
and its unit ball, ${\mathcal{ B}}_1 = {\left\{ \Phi\in{\mathcal{ M}}_{\infty} \; ; \; {\left\| \Phi \right\|_{{{\infty}}} }   \leq1 \right\} }$.
Now, ${{\mathcal{ M}}\left[{\mathbf{ X}};{\mathbb{C}}\right] }$ embeds into the dual space ${\mathcal{ M}}_{\infty}^*$ in a
natural way; \ endow it with the appropriate weak* topology.
The following results are straightforward:

\begin{lemma}{\sf \label{simple.dense}}  
  The {\bf simple functions} of the form $\Phi =  \sum_{n} \phi_n
{{{\mathsf{ 1\!\!1}}}_{{{\mathcal{ U}}_n}}}$ are dense in ${\mathcal{ M}}_{\infty}$
 (where $\phi_n\in{\mathbb{C}}$ and ${\mathcal{ U}}_n\subset{\mathbf{ X}}$ are
measurable).

The weak* topology on ${{\mathcal{ M}}\left[{\mathbf{ X}};{\mathbb{C}}\right] }$  is determined
by convergence on measurable sets.  Thus, a sequence
$\{\mu_n\}_{n=1}^{\infty}$ converges to $\mu\in{{\mathcal{ M}}\left[{\mathbf{ X}};{\mathbb{C}}\right] }$
if and only if $\displaystyle \mu_n[{\mathcal{ U}}] { -\!\!\!-\!\!\!-\!\!\!-\!\!\!\!\!\!\!\!\!\!\!  ^{{\scriptscriptstyle }}_{{\scriptscriptstyle n{\rightarrow}{\infty}}}   \!\!\!\!\!\!\!\!\!\longrightarrow }\mu[{\mathcal{ U}}]$ for all
measurable ${\mathcal{ U}}\subset{\mathbf{ X}}$.\hrulefill\ensuremath{\Box}
 \end{lemma}

\medskip

If a function ${\mathbf{ X}}\ni x\mapsto \mu_x\in{{\mathcal{ M}}\left[{\mathbf{ X}};{\mathbb{C}}\right] }$ is
measurable relative to the weak* Borel algebra of ${{\mathcal{ M}}\left[{\mathbf{ X}};{\mathbb{C}}\right] }$,
and $\nu$ is some other measure on ${\mathbf{ X}}$, then $\mu = \int_{\mathbf{ X}} \mu_x \
d\nu[x]$ is the measure so that, for all ${\mathcal{ U}}\subset{\mathbf{ X}}$, \ $\mu[{\mathcal{ U}}]
= \int_{\mathbf{ X}} \mu_x[{\mathcal{ U}}] \ d\nu[x]$; by Lemma \ref{simple.dense}, this
well-defines the action of $\mu$ on ${\mathcal{ M}}_{\infty}$.

  If ${\mathbf{ Q}}:{{\mathcal{ M}}\left[{\mathbf{ X}};{\mathbb{C}}\right] }{\,\raisebox{0.3em}{$-$}\!\!\!\!\!\!\raisebox{-0.3em}{$\leftarrow$}\!\!\!\!\supset}$, then define $\displaystyle{\left\| {\mathbf{ Q}} \right\|_{{}} }    :=
\sup_{x\in{\mathbf{ X}}}{\left\| {\mathbf{ q}}_x \right\|_{{var}} }   $ (note:  this is {\em not} the
operator norm of ${\mathbf{ Q}}$).  Say
that ${\mathbf{ Q}}$ is {\bf smooth} if ${\mathbf{ Q}}$ is linear, measurable relative to
the weak* Borel sigma algebra, and ${\left\| {\mathbf{ Q}} \right\|_{{}} }    < {\infty}$.

\begin{lemma}{\sf \label{smooth}}   \setcounter{enumi}{\thethm} \begin{list}{{\bf (\alph{enumii})}}{\usecounter{enumii}} 			{\setlength{\leftmargin}{0em} 			\setlength{\rightmargin}{0em}}
\item  If ${\mathbf{ Q}}$ is smooth, then its adjoint $\,^{\dagger}\!{{\mathbf{ Q}}}:{\mathcal{ M}}_{\infty}{\,\raisebox{0.3em}{$-$}\!\!\!\!\!\!\raisebox{-0.3em}{$\leftarrow$}\!\!\!\!\supset}$
is a well-defined, bounded linear operator, and ${\left\| \,^{\dagger}\!{{\mathbf{ Q}}} \right\|_{{{\infty}}} }   
\leq {\left\| {\mathbf{ Q}} \right\|_{{}} }   $.

\item If ${\mathbf{ X}}\ni x\mapsto {\mathbf{ q}}_x\in{{\mathcal{ M}}\left[{\mathbf{ X}};{\mathbb{C}}\right] }$ is a measurable
function and $M = \displaystyle\sup_{x\in{\mathbf{ X}}}{\left\| {\mathbf{ q}}_x \right\|_{{var}} }   <{\infty}$, then
the function ${\mathbf{ Q}}:{{\mathcal{ M}}\left[{\mathbf{ X}};{\mathbb{C}}\right] }{\,\raisebox{0.3em}{$-$}\!\!\!\!\!\!\raisebox{-0.3em}{$\leftarrow$}\!\!\!\!\supset}$ defined:
${\mathbf{ Q}}(\mu) = \int_{\mathbf{ X}} {\mathbf{ q}}_x \ d\mu[x]$ is smooth and
continuous, and ${\left\| {\mathbf{ Q}} \right\|_{{}} }   =M$.
\end{list}
 \end{lemma}
\bprf[Proof of (a):] For any $\phi\in{\mathcal{ M}}_{\infty}$, and any $x\in{\mathbf{ X}}$, define $(\,^{\dagger}\!{{\mathbf{ Q}}}\phi)(x)
\ = {\left\langle \phi,\ {\mathbf{ q}}_x \right\rangle }$.  Then $\,^{\dagger}\!{{\mathbf{ Q}}}(\phi)$ is measurable
(the function ${\mathbf{ X}}\ni x\mapsto \delta_x\in{{\mathcal{ M}}\left[{\mathbf{ X}};{\mathbb{C}}\right] }$ is
measurable; hence, so is the function $(x\mapsto
{\mathbf{ q}}_x)$; thus, so is $\,^{\dagger}\!{{\mathbf{ Q}}}(\phi)$).  Also, $\displaystyle
{\left\| \,^{\dagger}\!{{\mathbf{ Q}}}(\phi) \right\|_{{{\infty}}} }    \ \leq \ {\left\| \phi \right\|_{{{\infty}}} }   \cdot
\sup_{x\in{\mathbf{ X}}}{\left\| {\mathbf{ q}}_x \right\|_{{var}} }   $.

{\bf \hspace{-1em}  Proof of (b): \ \ } Clearly, ${\mathbf{ Q}}$ is  well-defined and linear,
and ${\left\| {\mathbf{ Q}} \right\|_{{}} }   =M$.    To see that ${\mathbf{ Q}}$ is continuous, let
${\mathcal{ U}}\subset{\mathbf{ X}}$; \ then ${\mathbf{ Q}}(\mu)[{\mathcal{ U}}] = \int_{\mathbf{ X}} {\mathbf{ q}}_x[{\mathcal{ U}}] \
d\mu[x]$.  The function ${\mathbf{ X}}\ni x\mapsto {\mathbf{ q}}_x[{\mathcal{ U}}]\in{\mathbb{C}}$ is
measurable; thus, if $\displaystyle\mu_n{ -\!\!\!-\!\!\!-\!\!\!-\!\!\!\!\!\!\!\!\!\!\!  ^{{\scriptscriptstyle }}_{{\scriptscriptstyle n{\rightarrow}{\infty}}}   \!\!\!\!\!\!\!\!\!\longrightarrow } \mu$ in the weak*
topology, then ${\mathbf{ Q}}(\mu_n)[{\mathcal{ U}}] \ = \ \int_{\mathbf{ X}} {\mathbf{ q}}_x[{\mathcal{ U}}] \ d\mu_n[x] \
{ -\!\!\!-\!\!\!-\!\!\!-\!\!\!\!\!\!\!\!\!\!\!  ^{{\scriptscriptstyle }}_{{\scriptscriptstyle n{\rightarrow}{\infty}}}   \!\!\!\!\!\!\!\!\!\longrightarrow } \int_{\mathbf{ X}} {\mathbf{ q}}_x[{\mathcal{ U}}] \ d\mu[x] \ = \ {\mathbf{ Q}}(\mu)[{\mathcal{ U}}]$.
 {\tt \hrulefill $\Box$ } \end{list}  \medskip  

  We define a {\bf Markov operator} to be a smooth linear operator
on ${{\mathcal{ M}}\left[{\mathbf{ X}};{\mathbb{C}}\right] }$ that maps ${{\mathcal{ M}}\left[{\mathbf{ X}}\right] }$ into itself.  By
Lemma \ref{smooth}, it suffices to define a measurable collection
$(x\mapsto{\mathbf{ q}}_x)$ of transition probability measures.  We will
be concerned with the following case:

        \refstepcounter{thm}        {\em }        \begin{list}{} 	{\setlength{\leftmargin}{1em} 	\setlength{\rightmargin}{0em}}         \item         {\bf Example  \thethm:} 
 Suppose ${\widetilde{\mathbb{M}}}={\mathbb{M}}\times{\mathbb{Z}}$ and
$\mu\in{{\mathcal{ M}}\left[{\mathcal{ A}}^{{\widetilde{\mathbb{M}}}}\right] }$ is an MRF
and consider the {\bf lamination process};  we claim the transition
probabilities are determined by a sequence $\{{\mathbf{ Q}}^{(n)}\}_{n\in{\mathbb{Z}}}$
of Markov operators.

 Let ${\mathbb{M}}_{k} := {\mathbb{M}}\times\{k\}$ for $k=n$ or $n+1$.  If ${\mathbf{ c}}\in{\mathcal{ A}}^{\widetilde{\mathbb{M}}}$ is a
$\mu$-random configuration, then for each ${\mathbf{ a}}\in{\mathcal{ A}}^{{\mathbb{M}}}$, \
${\mathbf{ q}}^{(n)}_{{\mathbf{ a}}}$ is the conditional distribution of
${\mathbf{ c}}\raisebox{-0.3em}{$\left|_{}\right.$}_{{\mathbb{M}}_{(n+1)}}$ given that ${\mathbf{ c}}\raisebox{-0.3em}{$\left|_{{\mathbb{M}}_n}\right.$} \ = \ {\mathbf{ a}}$.
Formally, if $\mu_{\mathbf{ a}}\in {{\mathcal{ M}}\left[{\mathcal{ A}}^{{\mathbb{M}}}\right] }$ is the conditional
distribution given ${\mathbf{ a}}$, then ${\mathbf{ q}}^{(n)}_{\mathbf{ a}} \ = \ {\mathbf{ pr}_{{{\mathbb{M}}_{(n+1)}}}}^*
\left(\mu_{\mathbf{ a}}\right)$.  The map ${\mathcal{ A}}^{\mathbb{M}}\ni
{\mathbf{ a}}\mapsto{\mathbf{ q}}^{(n)}_{\mathbf{ a}}\in{{\mathcal{ M}}\left[{\mathcal{ A}}^{\mathbb{M}}\right] }$ is measurable because the
map ${\mathcal{ A}}^{{\mathbb{M}}_n} \ni{\mathbf{ a}}\mapsto\mu_{\mathbf{ a}}\in
{{\mathcal{ M}}\left[{\mathcal{ A}}^{{\mathbb{M}}}\right] }$ is measurable \cite{SchwartzDisintegrate}, while 
${\mathbf{ pr}_{{{\mathbb{M}}_{(n+1)}}}}:{{\mathcal{ M}}\left[{\mathcal{ A}}^{{\mathbb{M}}}\right] }{{\longrightarrow}}{{\mathcal{ M}}\left[{\mathcal{ A}}^{{\mathbb{M}}_{(n+1)}}\right] }$ is
continuous.
  \hrulefill			        \end{list}   			 

  If ${\boldsymbol{\chi }}\in{\mathcal{ M}}_{\infty}$, then let ${{\boldsymbol{\chi }}}_\bullet:{\mathcal{ M}}_{\infty}{\,\raisebox{0.3em}{$-$}\!\!\!\!\!\!\raisebox{-0.3em}{$\leftarrow$}\!\!\!\!\supset}$ be the
bounded linear operator induced by multiplication with ${\boldsymbol{\chi }}$: for
any $\phi\in{\mathcal{ M}}_{\infty}$ and $x\in{\mathbf{ X}}$, \ $\left({{\boldsymbol{\chi }}}_\bullet\phi\right)(x) =
{\boldsymbol{\chi }}(x)\cdot\phi(x)$.  To establish that Markov processes on ${\mathcal{ A}}^{\mathbb{Z}}$ were EHM (Proposition \ref{markov.proc.exponential.HM}), we bounded the norm of
operators of the form ${{\boldsymbol{\xi }}}_\bullet\circ{\mathbf{ Q}}\circ{{\boldsymbol{\chi }}}_\bullet\circ{\mathbf{ P}}$, where
${\boldsymbol{\xi }},{\boldsymbol{\chi }}\in\widehat{{\mathcal{ A}}^{\mathbb{M}}}$.  We will employ a similar strategy to show that
MRFs are EHM; \ this will require the following result:

\begin{lemma}{\sf \label{matrix.for.triplex}}  
 Let ${\mathbf{ Q}},{\mathbf{ P}}:{{\mathcal{ M}}\left[{\mathbf{ X}}\right] }{\,\raisebox{0.3em}{$-$}\!\!\!\!\!\!\raisebox{-0.3em}{$\leftarrow$}\!\!\!\!\supset}$ be Markov operators.  For any $\phi\in{\mathcal{ M}}_{\infty}$, define $\mu^\phi_{x}\in{{\mathcal{ M}}\left[{\mathbf{ X}};\ {\mathbb{C}}\right] }$ \ by: \
$d\mu^\phi_{x} \ = \ \,^{\dagger}\!{{\mathbf{ P}}}(\phi)\cdot d{\mathbf{ q}}_x$.
 \setcounter{enumi}{\thethm} \begin{list}{{\bf (\alph{enumii})}}{\usecounter{enumii}} 			{\setlength{\leftmargin}{0em} 			\setlength{\rightmargin}{0em}}
 \item For any ${\boldsymbol{\chi }}\in{\mathcal{ M}}_{\infty}$, \ $\displaystyle
   {\left\| \,^{\dagger}\!{{\mathbf{ Q}}}\circ{{\boldsymbol{\chi }}}_\bullet\circ\,^{\dagger}\!{{\mathbf{ P}}} \right\|_{{{\infty}}} }   
\ = \ \sup_{\phi\in{\mathcal{ B}}_1} \sup_{x\in{\mathbf{ X}}} \left|{\left\langle {\boldsymbol{\chi }}, \ \mu^\phi_{x} \right\rangle }\right|$.

 \item Suppose that
$\mu\in{{\mathcal{ M}}\left[{\mathbf{ X}}^{{\mathbb{Z}}}\right] }$ is a Markov process and
${\mathbf{ Q}}$ and ${\mathbf{ P}}$ are the transition probability operators at time $0$ and
$1$, respectively.  For any $u,w\in{\mathbf{ X}}$, \ 
let $\mu_u = {\mathbf{ P}}\circ{\mathbf{ Q}}(\delta_u)$ be the conditional probability measure
on ${\mathbf{ X}}$ at time $2$ induced by state $u\in{\mathbf{ X}}$ at time $0$, and let
$\mu^{w}_{u}$ be the sandwich measure on ${\mathbf{ X}}$ induced by $u\in{\mathbf{ X}}$
at time $0$ and $w\in{\mathbf{ X}}$ at time $2$.  Then for any $\Phi\in{\mathcal{ M}}_{\infty}$,
\quad $\displaystyle \mu^\Phi_{u} \ = \ 
\int_{{\mathbf{ X}}} \Phi(w) \cdot \mu^{w}_{u}  { \;\; d\mu}_u\left[w\right]$.
\end{list}
 \end{lemma}
\bprf[Proof of (a):]
For any $\phi\in{\mathcal{ B}}_1$ and $x\in{\mathbf{ X}}$, \ 
$ \,^{\dagger}\!{{\mathbf{ Q}}}\circ{{\boldsymbol{\chi }}}_\bullet\circ\,^{\dagger}\!{{\mathbf{ P}}}(\phi)(x) \ = \ 
{\left\langle \,^{\dagger}\!{{\mathbf{ Q}}}\circ{{\boldsymbol{\chi }}}_\bullet\circ\,^{\dagger}\!{{\mathbf{ P}}}(\phi), \ \delta_x \right\rangle }
\ = \ 
{\left\langle {{\boldsymbol{\chi }}}_\bullet\circ\,^{\dagger}\!{{\mathbf{ P}}}(\phi), \ {\mathbf{ Q}}(\delta_x) \right\rangle }
\ = \ 
{\left\langle {\boldsymbol{\chi }} \cdot \,^{\dagger}\!{{\mathbf{ P}}}(\phi), \ {\mathbf{ q}}_x \right\rangle }
\ = \ 
{\left\langle {\boldsymbol{\chi }}, \ \mu^\phi_{x} \right\rangle }$.  \  Thus, 
\begin{eqnarray*}
 {\left\| \,^{\dagger}\!{{\mathbf{ Q}}}\circ{{\boldsymbol{\chi }}}_\bullet\circ\,^{\dagger}\!{{\mathbf{ P}}} \right\|_{{{\infty}}} }   
& = & \sup_{\phi\in{\mathcal{ B}}_1} {\left\| \,^{\dagger}\!{{\mathbf{ Q}}}\circ{{\boldsymbol{\chi }}}_\bullet\circ\,^{\dagger}\!{{\mathbf{ P}}}(\phi) \right\|_{{{\infty}}} }   
\quad = \quad  \sup_{\phi\in{\mathcal{ B}}_1} \sup_{x\in{\mathbf{ X}}} \left| \,^{\dagger}\!{{\mathbf{ Q}}}\circ{{\boldsymbol{\chi }}}_\bullet\circ\,^{\dagger}\!{{\mathbf{ P}}}(\phi)(x)\right|
\\ &=&
 \sup_{\phi\in{\mathcal{ B}}_1} \sup_{x\in{\mathbf{ X}}} \left| {\left\langle {\boldsymbol{\chi }}, \ \mu^\phi_{x} \right\rangle } \right|.
\end{eqnarray*}
{\bf \hspace{-1em}  Proof of (b): \ \ }  We want to show that
 for any ${\mathbf{ V}}\subset{\mathbf{ X}}$, 
\[
 \mu^\Phi_{u}({\mathbf{ V}}) \quad = \quad 
\int_{{\mathbf{ X}}} \Phi(w) \cdot \mu^{w}_{u}({\mathbf{ V}})  { \;\; d\mu}_u\left[w\right].
\]
First, suppose that ${\mathbf{ W}}\subset{\mathbf{ X}}$
and $\Phi={{{\mathsf{ 1\!\!1}}}_{{{\mathbf{ W}}}}}$.  Let $\mu^{\mathbf{ W}}_x := \mu^{{{{\mathsf{ 1\!\!1}}}_{{{\mathbf{ W}}}}}}_x$; \ 
thus, $d\mu^{\mathbf{ W}}_x = \ \,^{\dagger}\!{{\mathbf{ p}}}\rule[-0.5em]{0em}{1em}_{\mathbf{ W}} \, d{\mathbf{ q}}_x$. Then:
\begin{eqnarray*}
\displaystyle \mu^{{\mathbf{ W}}}_{u}({\mathbf{ V}})
& = &  \int_{{\mathbf{ V}}}  d\mu^{{\mathbf{ W}}}_{u}[v]
 \quad =  \quad \int_{{\mathbf{ V}}} \,^{\dagger}\!{{\mathbf{ p}}}\rule[-0.5em]{0em}{1em}_{{\mathbf{ W}}}\left(v\right) \ d{\mathbf{ q}}_u\left[v\right]
\quad=\quad
 \int_{{\mathbf{ V}}} {\mathbf{ p}}_{v}\left[{\mathbf{ W}}\right] \ d{\mathbf{ q}}_u\left[v\right]
\\ &=&
 \int_{{\mathbf{ V}}} \mu\left[\frac{x_{2} \in {\mathbf{ W}}}{  x_1 = v}\right] \ d{\mathbf{ q}}_{u}\left[v\right]
\quad\raisebox{-1ex}{$\overline{\overline{{\scriptscriptstyle{\mathrm{(a)}}}}}$}\quad
 \int_{{\mathbf{ V}}} \mu\left[\frac{x_2\in {\mathbf{ W}}}{ (x_1= v) \&(x_0=u)} \right]  \ d{\mathbf{ q}}_{u}\left[v\right] 
\\ & \raisebox{-1ex}{$\overline{\overline{{\scriptscriptstyle{\mathrm{(b)}}}}}$} &
 \int_{\mathbf{ W}}  \mu_u^w({\mathbf{ V}}) \ d\mu_u[w]
\quad=\quad
\int_{\mathbf{ X}} {{{\mathsf{ 1\!\!1}}}_{{{\mathbf{ W}}}}}(w)  \cdot  \mu_u^w({\mathbf{ V}}) \ d\mu_u[w].
\end{eqnarray*}
  {\bf(a)} follows from the Markov property.
  To see {\bf(b)}, let ${\mathcal{ X}}_k$ be the sigma-subalgebra of ${\mathbf{ X}}^{\mathbb{Z}}$
generated by coordinate $x_k$, and let $\ensuremath{\mathbf{ E}_{k}\left[ \bullet \right]}$ (resp. $\ensuremath{\mathbf{ E}_{k,j}\left[ \bullet \right]}$) be the
conditional expectation with respect to ${\mathcal{ X}}_k$ (resp. ${\mathcal{ X}}_k\vee{\mathcal{ X}}_j$).
Let
${\mathbf{ W}}_2={\left\{ {\mathbf{ x}}\in{\mathbf{ X}}^{\mathbb{Z}} \; ; \; x_2\in{\mathbf{ W}} \right\} }$ and
${\mathbf{ V}}_1={\left\{ {\mathbf{ x}}\in{\mathbf{ X}}^{\mathbb{Z}} \; ; \; x_1\in{\mathbf{ V}} \right\} }$.  Then 
\begin{eqnarray*} 
\lefteqn{\int_{{\mathbf{ V}}}
\mu\left[\frac{x_2\in {\mathbf{ W}}}{ (x_1= v) \&(x_0=u)} \right] \ d{\mathbf{ q}}_{u}\left[v\right]
\quad= \quad\int_{{\mathbf{ V}}} \ensuremath{\mathbf{ E}_{0,1}\left[ {{{\mathsf{ 1\!\!1}}}_{{{\mathbf{ W}}_2}}} \right]}(u,v) \ d{\mathbf{ q}}_{u}\left[v\right]} \\
& = & \int_{\mathbf{ X}} {{{\mathsf{ 1\!\!1}}}_{{{\mathbf{ V}}}}}\cdot\ensuremath{\mathbf{ E}_{0,1}\left[ {{{\mathsf{ 1\!\!1}}}_{{{\mathbf{ W}}_2}}} \right]}(u,v) \ d{\mathbf{ q}}_{u}\left[v\right]
\quad= \quad\ensuremath{\mathbf{ E}_{0}\left[ \rule[-0.5em]{0em}{1em} {{{\mathsf{ 1\!\!1}}}_{{{\mathbf{ V}}_1}}}\cdot\ensuremath{\mathbf{ E}_{0,1}\left[ {{{\mathsf{ 1\!\!1}}}_{{{\mathbf{ W}}_2}}} \right]} \right]}(u)\\
& = & \ensuremath{\mathbf{ E}_{0}\left[ \rule[-0.5em]{0em}{1em}\ensuremath{\mathbf{ E}_{0,1}\left[ {{{\mathsf{ 1\!\!1}}}_{{{\mathbf{ V}}_1}}}\cdot {{{\mathsf{ 1\!\!1}}}_{{{\mathbf{ W}}_2}}} \right]} \right]}(u)
\quad= \quad\ensuremath{\mathbf{ E}_{0}\left[ {{{\mathsf{ 1\!\!1}}}_{{{\mathbf{ V}}_1}}}\cdot {{{\mathsf{ 1\!\!1}}}_{{{\mathbf{ W}}_2}}} \right]}(u) \\
& = & \ensuremath{\mathbf{ E}_{0}\left[ \rule[-0.5em]{0em}{1em}\ensuremath{\mathbf{ E}_{0,2}\left[ {{{\mathsf{ 1\!\!1}}}_{{{\mathbf{ V}}_1}}}\cdot {{{\mathsf{ 1\!\!1}}}_{{{\mathbf{ W}}_2}}} \right]} \right]}(u)
\quad = \quad\ensuremath{\mathbf{ E}_{0}\left[ \rule[-0.5em]{0em}{1em} {{{\mathsf{ 1\!\!1}}}_{{{\mathbf{ W}}_2}}}\cdot\ensuremath{\mathbf{ E}_{0,2}\left[ {{{\mathsf{ 1\!\!1}}}_{{{\mathbf{ V}}_1}}} \right]} \right]}(u)\\
& = & \int_{{\mathbf{ W}}} \mu\left[\frac{x_1\in {\mathbf{ V}}}{ (x_2= w) \&(x_0=u)} \right]  \ d\mu_{u}\left[w\right]
\ = \ \int_{\mathbf{ W}}  \mu_u^w({\mathbf{ V}}) \ d\mu_u[w].\end{eqnarray*}

  Next, if $\displaystyle\Phi = \sum_n \phi_n {{{\mathsf{ 1\!\!1}}}_{{{\mathbf{ W}}_n}}}$ is a simple function,
then $\displaystyle\,^{\dagger}\!{{\mathbf{ P}}}(\Phi) =  \sum_n \phi_n \,^{\dagger}\!{{\mathbf{ p}}}\rule[-0.5em]{0em}{1em}_{{\mathbf{ W}}_n}$,
so that $\displaystyle d\mu^\Phi_{x} \ = \ \,^{\dagger}\!{{\mathbf{ P}}}(\Phi)\cdot\  d{\mathbf{ q}}_x
\ =  \  \sum_n \phi_n \,^{\dagger}\!{{\mathbf{ p}}}\rule[-0.5em]{0em}{1em}_{{\mathbf{ W}}_n} \ d{\mathbf{ q}}_x \ = \ 
 \sum_n \phi_n \ d\mu^{{\mathbf{ W}}_n}_{x}$.
Thus, 
\begin{eqnarray*}
\mu^\Phi_{x}({\mathbf{ V}})
&=&  \int_{\mathbf{ V}} d\mu^\Phi_{x}
\quad = \quad  \sum_n \phi_n \cdot \int_{\mathbf{ V}} d\mu^{{\mathbf{ W}}_n}_{x}
\quad = \quad  \sum_n \phi_n \cdot \mu^{{\mathbf{ W}}_n}_{x}[{\mathbf{ V}}]
\\ &=& 
   \sum_n \phi_n \cdot  \int_{{\mathbf{ X}}}{{{\mathsf{ 1\!\!1}}}_{{{\mathbf{ W}}_n}}}(w) \mu^{w}_{u}({\mathbf{ V}}) 
{ \;\; d\mu}_u\left[w\right] 
\quad = \quad
   \int_{\mathbf{ X}}  \left(\sum_n \phi_n {{{\mathsf{ 1\!\!1}}}_{{{\mathbf{ W}}_n}}}(w) \right)\cdot \mu^{w}_{u}({\mathbf{ V}})
  { \;\; d\mu}_u\left[w\right]
\\ &=&
 \int_{\mathbf{ X}} \Phi(w) \cdot \mu^{w}_{u}({\mathbf{ V}})  { \;\; d\mu}_u\left[w\right],
\qquad\mbox{as desired.}
\end{eqnarray*}
Finally, if $\{\Phi_n\}_{n=1}^{\infty}$ is a bounded sequence of simple functions
so that  $\displaystyle\Phi_n { -\!\!\!-\!\!\!-\!\!\!-\!\!\!\!\!\!\!\!\!\!\!  ^{{\scriptscriptstyle }}_{{\scriptscriptstyle n{\rightarrow}{\infty}}}   \!\!\!\!\!\!\!\!\!\longrightarrow }\Phi $ in ${\mathcal{ M}}_{\infty}$, then
$\displaystyle\,^{\dagger}\!{{\mathbf{ P}}}(\Phi_n) { -\!\!\!-\!\!\!-\!\!\!-\!\!\!\!\!\!\!\!\!\!\!  ^{{\scriptscriptstyle }}_{{\scriptscriptstyle n{\rightarrow}{\infty}}}   \!\!\!\!\!\!\!\!\!\longrightarrow } \,^{\dagger}\!{{\mathbf{ P}}}(\Phi)$ in ${\mathcal{ M}}_{\infty}$, 
so that $\displaystyle\mu^{\Phi_n}_{x} { -\!\!\!-\!\!\!-\!\!\!-\!\!\!\!\!\!\!\!\!\!\!  ^{{\scriptscriptstyle }}_{{\scriptscriptstyle n{\rightarrow}{\infty}}}   \!\!\!\!\!\!\!\!\!\longrightarrow } \mu^\Phi_{x}$ 
in the weak* topology on ${{\mathcal{ M}}\left[{\mathbf{ X}}\right] }$.  But by dominated convergence,
we also know that
\[
 \mu^{\Phi_n}_{x}({\mathbf{ V}}) 
\quad = \quad 
\int_{\mathbf{ X}} \Phi_n(w) \cdot \mu^{w}_{u}({\mathbf{ V}})  { \;\; d\mu}_u\left[w\right] 
\quad { -\!\!\!-\!\!\!-\!\!\!-\!\!\!\!\!\!\!\!\!\!\!  ^{{\scriptscriptstyle }}_{{\scriptscriptstyle n{\rightarrow}{\infty}}}   \!\!\!\!\!\!\!\!\!\longrightarrow }\quad 
\int_{\mathbf{ X}} \Phi(w) \cdot \mu^{w}_{u}({\mathbf{ V}})  { \;\; d\mu}_u\left[w\right]; 
\] 
hence, $\displaystyle\mu^\Phi_{x}({\mathbf{ V}}) \ = \ \int_{\mathbf{ X}} \Phi(w) \cdot
\mu^{w}_{u}({\mathbf{ V}}) { \;\; d\mu}_u\left[w\right]$, for all ${\mathbf{ V}}\subset{\mathbf{ X}}$, as
desired.
 {\tt \hrulefill $\Box$ } \end{list}  \medskip  
\subsection{Uniform Harmonic Mixing
\label{S:uhm}}

  Now, let ${\mathbf{ X}} = {\mathcal{ A}}^{\mathbb{M}}$, and consider a $1$-step Markov process
on ${\mathcal{ A}}^{\mathbb{M}}$ determined by a sequence of Markov operators
${\left\{ {\mathbf{ Q}}^{(n)} \; ; \; n\in{\mathbb{Z}} \right\} }$.  For every $n\in{\mathbb{Z}}$, ${\mathbf{ a}}\in{\mathcal{ A}}^{\mathbb{M}}$ and
 $\phi\in{\mathcal{ M}}_{\infty}$, define $\,_{n\!}\mu^\phi_{{\mathbf{ a}}}\in{{\mathcal{ M}}\left[{\mathcal{ A}}^{\mathbb{M}}; \ {\mathbb{C}}\right] }$ so that
$d\,_{n\!}\mu^\phi_{{\mathbf{ a}}} \ = \
\,^{\dagger}\!{{\mathbf{ Q}}}^{(n+1)}(\phi) \, d{\mathbf{ q}}^{(n)}_{\mathbf{ a}}$ as in Lemma
\ref{matrix.for.triplex}.

If $\lambda>0$  then
the sequence
${\left\{ {\mathbf{ Q}}^{(n)} \; ; \; n\in{\mathbb{Z}} \right\} }$ is {\bf uniformly harmonically mixing}
with {\bf decay parameter} $\lambda$ (or ``$\lambda$-UHM'') if, for every
${\mathbf{ a}}\in{\mathcal{ A}}^{\mathbb{M}}$ and measurable $\phi\in{\mathcal{ B}}_1$, and every $n\in{\mathbb{Z}}$, the
measure $\,_{n\!}\mu^\phi_{{\mathbf{ a}}}$ is $\lambda$-EHM.  Thus, applying
Lemma \ref{matrix.for.triplex}(a), we have:
$\displaystyle {\left\|  \,^{\dagger}\!{{\mathbf{ Q}}}^{(n+1)} \circ {{\boldsymbol{\chi }}}_\bullet \circ \,^{\dagger}\!{{\mathbf{ Q}}}^{(n)} \right\|_{{{\infty}}} }   
\ \leq \ e^{-\lambda\cdot R}$ for any ${\boldsymbol{\chi }}\in\widehat{{\mathcal{ A}}^{\mathbb{M}}}$ 
with ${{\sf rank}\left[{\boldsymbol{\chi }}\right]}\geq R$.

\begin{prop}{\sf \label{sandwich.unto.UHM}}  
 Let ${\widetilde{\mathbb{U}}}={\mathbb{U}}\times\{-1,0,1\}$ and $\mu\in{{\mathcal{ M}}\left[{\mathcal{ A}}^{{\mathbb{M}}\times{\mathbb{Z}}}\right] }$ be a ${\widetilde{\mathbb{U}}}$-MRF such that  
$\forall n\in{\mathbb{Z}}$, \ 
 ${\mathbf{ a}}\in{\mathcal{ A}}^{{\mathbb{M}}\times\{n\}}$,  and
${\mathbf{ c}}\in{\mathcal{ A}}^{{\mathbb{M}}\times\{n+2\}}$, the sandwich measure $\mu_{\mathbf{ a}}^{\mathbf{ c}}$ is
$\lambda$-EHM.  Then  ${\left\{ {\mathbf{ Q}}^{\left(n\right)} \; ; \; n\in{\mathbb{Z}} \right\} }$ is
$\lambda$-UHM.   \end{prop}
\bprf Let $\phi\in{\mathcal{ B}}_1$.  By Lemma \ref{matrix.for.triplex}(b), \ 
$\displaystyle \,_{n\!}\mu^\phi_{{\mathbf{ a}}} \ = \ 
\int_{{\mathcal{ A}}^{\mathbb{M}}} \phi({\mathbf{ c}})\cdot  \mu^{{\mathbf{ c}}}_{{\mathbf{ a}}}  { \;\; d\mu}_{\mathbf{ a}}\left[{\mathbf{ c}}\right]$,
where $\mu_{\mathbf{ a}} = {\mathbf{ Q}}^{(n+1)}\circ{\mathbf{ Q}}^{(n)}(\delta_{\mathbf{ a}})$. \ 
By hypothesis, $\mu_{\mathbf{ a}}^{\mathbf{ c}}$ is $\lambda$-EHM for all ${\mathbf{ c}}\in{\mathcal{ A}}^{\mathbb{M}}$; \ 
apply Lemma \ref{integral.ehm} to conclude that $\,_{n\!}\mu^\phi_{{\mathbf{ a}}}$
is also $\lambda$-EHM.
 {\tt \hrulefill $\Box$ } \end{list}  \medskip

\begin{prop}{\sf \label{uniform.harmonic.mixing}}  
 Let $\mu\in{{\mathcal{ M}}\left[{\mathcal{ A}}^{{\mathbb{M}}\times{\mathbb{Z}}}\right] }$ be a the path distribution of
a Markov process determined by
Markov operators ${\left\{ {\mathbf{ Q}}^{(n)} \; ; \; n\in{\mathbb{Z}} \right\} }$.

 If ${\left\{ {\mathbf{ Q}}^{(n)} \; ; \; n\in{\mathbb{Z}} \right\} }$ is $\lambda$-UHM, 
then $\mu$ is $\lambda'$-HM, where $\lambda'=\lambda/2$.
 \end{prop}
\bprf
Let ${\boldsymbol{\chi }} \in \widehat{{\mathcal{ A}}^{{\mathbb{M}}\times{\mathbb{Z}}}}$ with ${{\sf rank}\left[{\boldsymbol{\chi }}\right]} = 2R$,
and suppose that $\displaystyle {\boldsymbol{\chi }} = \bigotimes_{k=0}^{2K} {\boldsymbol{\chi }}^{(k)}$,
where, for all $k\in{\left[ 0..2K \right]}$,
${\boldsymbol{\chi }}^{(k)}$ is a character on ${\mathcal{ A}}^{{\mathbb{M}}\times\{n_k\}}$, 
with ${{\sf rank}\left[{\boldsymbol{\chi }}^{(k)}\right]} = R_k$, for some
$n_0 < n_1 < \ldots < n_{2K}$.
For any $k\in{\left[ 1..2K \right]}$, define $\,^{\dagger}\!{{\mathbf{ Q}}}_k =
\,^{\dagger}\!{{\mathbf{ Q}}}^{(n_k)} \circ \,^{\dagger}\!{{\mathbf{ Q}}}^{(n_k-1)} \circ 
\ldots \circ \,^{\dagger}\!{{\mathbf{ Q}}}^{(n_{(k-1)}+2)}
 \circ \,^{\dagger}\!{{\mathbf{ Q}}}^{(n_{(k-1)}+1)}$. 

If ${\left\{ \eta_n \; ; \; n\in{\mathbb{Z}} \right\} } \subset {{\mathcal{ M}}\left[{\mathcal{ A}}^{\mathbb{M}}\right] }$ are the 
state distributions of the process, then it is not hard
to show:
\begin{eqnarray*}
 {\left\langle {\boldsymbol{\chi }},\ \mu \right\rangle }
& = & 
\left\langle \rule[-0.5em]{0em}{1em}
\,^{\dagger}\!{{\mathbf{ Q}}}^{\left(n_{(2K)}+1\right)} \circ {{\boldsymbol{\chi }}_\bullet^{(2K)}} \circ  
\,^{\dagger}\!{{\mathbf{ Q}}}_{2K}  \circ {{\boldsymbol{\chi }}_\bullet^{(2K-1)}} \circ \,^{\dagger}\!{{\mathbf{ Q}}}_{(2K-1)}  \circ
\cdots \circ
\,^{\dagger}\!{{\mathbf{ Q}}}_2  \circ {{\boldsymbol{\chi }}_\bullet^{(1)}} \circ
\,^{\dagger}\!{{\mathbf{ Q}}}_1 \left({\boldsymbol{\chi }}^{(0)}\right), \right.
\\ &&\hspace{5em} 
\left.
\ \eta_{\left(n_{(2K)}+1\right)} \rule[-0.5em]{0em}{1em} \right\rangle,
\end{eqnarray*}
 (see e.g. {\bf Claim 1} of Proposition 8 in \cite{PivatoYassawi1}).
Thus,
\begin{eqnarray*}
 \lefteqn{\left|{\left\langle {\boldsymbol{\chi }},\ \mu \right\rangle }\right|}\\
& \raisebox{-1ex}{${{\leq} \atop {\scriptscriptstyle{\mathrm{(1)}}}}$} & 
{\left\| \,^{\dagger}\!{{\mathbf{ Q}}}^{\left(n_{(2K)}+1\right)} \circ {{\boldsymbol{\chi }}_\bullet^{(2K)}} \circ  
\,^{\dagger}\!{{\mathbf{ Q}}}_{2K}  \circ {{\boldsymbol{\chi }}_\bullet^{(2K-1)}} \circ
\ldots \circ
\,^{\dagger}\!{{\mathbf{ Q}}}_{2}  \circ {{\boldsymbol{\chi }}_\bullet^{(1)}} \circ
\,^{\dagger}\!{{\mathbf{ Q}}}_{1} \left({\boldsymbol{\chi }}^{(0)}\right) \right\|_{{{\infty}}} }   \\
& \raisebox{-1ex}{${{\leq} \atop {\scriptscriptstyle{\mathrm{(2)}}}}$} &
{\left\| \,^{\dagger}\!{{\mathbf{ Q}}}^{\left(n_{(2K)}+1\right)} \circ {{\boldsymbol{\chi }}_\bullet^{(2K)}} \circ  
\,^{\dagger}\!{{\mathbf{ Q}}}_{2K}  \circ {{\boldsymbol{\chi }}_\bullet^{(2K-1)}} \circ
\ldots \circ
\,^{\dagger}\!{{\mathbf{ Q}}}_{2}  \circ {{\boldsymbol{\chi }}_\bullet^{(1)}} \circ
\,^{\dagger}\!{{\mathbf{ Q}}}_{1} \right\|_{{{\infty}}} }   \\
& \raisebox{-1ex}{${{\leq} \atop {\scriptscriptstyle{\mathrm{(3)}}}}$} & 
{\left\| \,^{\dagger}\!{{\mathbf{ Q}}}^{\left(n_{(2K)}+1\right)} \circ {{\boldsymbol{\chi }}_\bullet^{(2K)}}
 \circ \,^{\dagger}\!{{\mathbf{ Q}}}_{2K} \right\|_{{{\infty}}} }    \cdot
\prod_{k=1}^{K-1} {\left\|  \,^{\dagger}\!{{\mathbf{ Q}}}_{(2k+1)}  \circ 
{{\boldsymbol{\chi }}_\bullet^{(2k)}} \circ \,^{\dagger}\!{{\mathbf{ Q}}}_{2k}   \right\|_{{{\infty}}} }   
\cdot \prod_{k=0}^{K-1} {\left\| {\boldsymbol{\chi }}_\bullet^{(2k+1)} \right\|_{{{\infty}}} }   \\
& \raisebox{-1ex}{${{\leq} \atop {\scriptscriptstyle{\mathrm{(4)}}}}$} & 
\prod_{k=1}^K {\left\| \,^{\dagger}\!{{\mathbf{ Q}}}^{\left(n_{(2k)}+1\right)}  \circ {{\boldsymbol{\chi }}_\bullet^{(2k)}} \circ
\,^{\dagger}\!{{\mathbf{ Q}}}^{\left(n_{(2k)}\right)} \right\|_{{{\infty}}} }    \\
& \raisebox{-1ex}{${{\leq} \atop {\scriptscriptstyle{\mathrm{(5)}}}}$} &
\prod_{k=1}^K \exp\left[-\lambda\cdot R_{2k}\right] \
\ = \ \exp\left[\sum_{k=1}^K -\lambda\cdot R_{2k} \right]
\ = \
 \exp\left[-\lambda\cdot \sum_{k=1}^K R_{2k} \right]
\end{eqnarray*}

{\bf(1)} Because $\eta_{\left(n_{(2K)}+1\right)}$ is
a probability measure.
\quad
{\bf(2)} Because ${\left\| {\boldsymbol{\chi }}^{(0)} \right\|_{{{\infty}}} }    = 1$. 

{\bf(3)} Separating out ${{\boldsymbol{\chi }}_\bullet^{(k)}}$ for all odd $k$.

{\bf(4)} Dropping ${{\boldsymbol{\chi }}_\bullet^{(k)}}$ for all odd $k$, and
$\,^{\dagger}\!{{\mathbf{ Q}}}^{\left(n_{(2k)}-1\right)},\ \,^{\dagger}\!{{\mathbf{ Q}}}^{\left(n_{(2k)}-2\right)}, \ \ldots,
\,^{\dagger}\!{{\mathbf{ Q}}}^{\left(n_{(2k-1)}+2\right)}$ for every $k$ (because
${\left\| \,^{\dagger}\!{{\mathbf{ Q}}}^{\left(n\right)} \right\|_{{}} }   \leq1$ for every $n\in{\mathbb{Z}}$).

{\bf(5)}  By UHM hypothesis and Lemma \ref{matrix.for.triplex}{\bf(a)}.
By the same logic, \ $\displaystyle
\left|{\left\langle {\boldsymbol{\chi }},\ \mu \right\rangle }\right|
\ \leq \
 \exp\left[-\lambda\cdot \sum_{k=1}^K R_{(2k-1)} \right]$.

 Now, clearly, one of \ $\displaystyle\left( \sum_{k=1}^K R_{2k}\right)$ and \ $\displaystyle\left(
\sum_{k=1}^K R_{(2k-1)}\right)$ must equal or exceed $R$, since together,
they sum to ${{\sf rank}\left[{\boldsymbol{\chi }}\right]} = 2R$.  Thus either \ $\displaystyle
-\lambda\cdot \sum_{k=1}^K R_{2k} \leq -\lambda\cdot R$ \ or \ $\displaystyle
-\lambda\cdot \sum_{k=1}^K R_{(2k-1)} \leq -\lambda\cdot R$.  \ Hence \ $\displaystyle
\left|{\left\langle {\boldsymbol{\chi }},\ \mu \right\rangle }\right|
\ \leq \
 -\lambda\cdot R$. 
 {\tt \hrulefill $\Box$ } \end{list}  \medskip

\begin{list}{} 			{\setlength{\leftmargin}{1em} 			\setlength{\rightmargin}{0em}}                         \item {\bf \hspace{-1em}  Proof of Proposition \ref{sandwich.unto.mixing}: \ \ }
  If $\mu$ is a MRF and all sandwich measures of $\mu$
are $\lambda$-EHM, then, by Proposition
\ref{sandwich.unto.UHM}, the sequence
${\left\{ {\mathbf{ Q}}^{\left(n\right)} \; ; \; n\in{\mathbb{Z}} \right\} }$ is $\lambda$-UHM.  Then, by
Proposition \ref{uniform.harmonic.mixing}, $\mu$ is $\lambda'$-HM, where
$\lambda' = \lambda/2$.
\hrulefill\ensuremath{\Box}\end{list}

\section{Conclusion}

 We have demonstrated that a broad class of probability measures on
${\mathcal{ A}}^{\mathbb{M}}$ weak*-converge to Haar measure in density, when acted on by
a wide class of LCA.  Many interesting questions remain.  For example,
can we establish harmonic mixing for Markov random fields {\em
without} full support, such as those supported on subshifts of finite
type?  What about measures on sofic shifts?  How can we characterize
either diffusion or harmonic mixing when ${\mathbb{M}}$ is not a lattice, but
instead, a nonabelian group or monoid?  Finally, what happens when
${\mathcal{ A}}$ is a nonabelian group?  The natural analogy of LCA for
nonabelian ${\mathcal{ A}}$ are ``multiplicative'' cellular automata
\cite{PivatoNACA}, where the local map is computed by
(noncommutatively) multiplying the values of neighbouring coordinates.
What is the asymptotic behaviour of measures under such automata?
{\small
\bibliographystyle{plain}
\bibliography{bibliography}
}

\end{document}